\documentclass[10pt]{amsart}

\usepackage{amsmath, amssymb}
\usepackage{amsfonts}
\usepackage{mathrsfs}
\usepackage[color, arrow,matrix,curve,cmtip,ps]{xy}
\usepackage{paralist}
\usepackage{amsthm}
\usepackage{tikz-cd}
\usetikzlibrary{arrows,calc,matrix}
\tikzset{
curvarr/.style={
  to path={ -- ([xshift=2ex]\tikztostart.east)
    |- (#1) [near end]\tikztonodes
    -| ([xshift=-2ex]\tikztotarget.west)
    -- (\tikztotarget)}
  }
}
\tikzset{%
    symbol/.style={%
        draw=none,
        every to/.append style={%
            edge node={node [sloped, allow upside down, auto=false]{$#1$}}}
    }
}

\usepackage[final]{pdfpages}
\usepackage{dsfont,txfonts}
\usepackage{tikz}
\usepackage{rotating}
\usepackage{mathrsfs}
\PassOptionsToPackage{hyphens}{url}\usepackage{hyperref}
\usepackage{enumitem}
\usepackage{graphicx}
\usepackage{mathtools}
\usetikzlibrary{arrows,chains,matrix,positioning,scopes}

\usepackage[T1]{fontenc}
\usepackage[variablett]{lmodern}
\usepackage{xcolor}
\usepackage[citestyle=alphabetic,bibstyle=alphabetic]{biblatex}
\usepackage{lipsum}

\newtheorem{theorem}{Theorem}[section]
\theoremstyle{definition}
\newtheorem{lemma}[theorem]{Lemma}

\newtheorem{proposition}[theorem]{Proposition}
\newtheorem{corollary}[theorem]{Corollary}
\newtheorem{definition}[theorem]{Definition}

\newtheorem{remark}[theorem]{Remark}

\newtheorem{example}[theorem]{Example}

\everymath{\displaystyle}
\allowdisplaybreaks

\newcommand{\set}[1]{\left\lbrace #1 \right\rbrace}

\newcommand{\B}{\mathbf{B}}
\newcommand{\E}{\mathbf{E}}
\newcommand{\G}{\mathbf{G}}
\newcommand{\pt}{\mathrm{pt}}
\newcommand{\Map}{\mathrm{Map}}
\newcommand{\Fun}{\mathrm{Fun}}
\newcommand{\Spc}{\mathrm{Spc}}
\newcommand{\Set}{\mathrm{Set}}
\newcommand{\Fin}{\mathrm{Fin}}
\newcommand{\FinInj}{\mathrm{FinInj}}
\newcommand{\Grpd}{\mathrm{Grpd}}
\newcommand{\Cat}{\mathrm{Cat}}
\newcommand{\Catinf}{\mathrm{Cat}_\infty}
\newcommand{\PrL}{\mathrm{Pr}^L}

\newcommand{\id}{\mathrm{id}}
\newcommand{\colim}{\mathrm{colim}}
\newcommand{\cof}{\mathrm{cof}}
\newcommand{\EMon}[1]{\mathbb{E}_{#1}\mbox{-}\mathrm{Mon}}
\newcommand{\EGrp}[1]{\mathbb{E}_{#1}\mbox{-}\mathrm{Grp}}
\newcommand{\EAlg}[1]{\mathbb{E}_{#1}\mbox{-}\mathrm{Alg}}

\newcommand{\real}[1]{\lvert #1 \rvert}
\newcommand{\core}{\mathrm{core}}

\newcommand{\rotadj}{\rotatebox[origin=c]{-90}{$\dashv$}}
\newcommand{\gp}{\mathrm{gp}}
\newcommand{\op}{\mathrm{op}}

\begin{document}
\title{Group completion via the action $\infty$-category}
\author{Georg Lehner}
\begin{abstract}
We give a generalization of Quillen's $S^{-1}S$ construction for arbitrary $\mathbb{E}_n$-monoids as an $\mathbb{E}_{n-1}$-monoidal $\infty$-category and show that its realization models the group completion provided that $n \geq 2$. We will also show how this construction is related to a variety of other constructions of the group completion. 
\end{abstract}
\maketitle

\section{Introduction}

After Quillen's original construction of higher algebraic $K$-theory \cite{10.1007/BFb0067053}, many different constructions to group complete a variety of sufficiently commutative higher algebraic gadgets, such as a homotopy commutative topological monoid, or a symmetric monoidal groupoid, have been considered in the literature. These constructions have been of fundamental importance in the study of higher $K$-theory as well as $L$-theory. We want to highlight the following examples.

\begin{enumerate}
\item If $(M,+)$ is an ordinary commutative monoid, we can put an equivalence relation on the set $M \times M$, by defining
$$ (a,b) \sim (c,d) \iff \exists k \in M : a+d+k = c+b+k. $$
The resulting quotient $M \times M/\sim$ will be an abelian group, and it is the initial abelian group under $M$, also called the Grothendieck group of $M$.
\item Suppose $(S,\oplus)$ is symmetric monoidal groupoid. Quillen \cite{10.1007/BFb0080003} constructs a symmetric monoidal category $S^{-1}S$ with objects given by pairs $(x,y) \in S$, and morphisms being equivalence classes of
$$(x_1,y_1) \xrightarrow{k,\alpha,\beta} (x_2,y_2)$$
consisting of $k \in S$ together with isomorphisms $\alpha : k \oplus x_1 \cong x_2$ and $\beta : k \oplus y_1 \cong y_2$. Two morphisms $(k,\alpha,\beta)$ and $(k',\alpha',\beta')$ are considered equivalent if there exists an isomorphism $\gamma : k \cong k'$ in $S$ that commutes with the structure involved. Quillen shows that in the case that $S = \text{Proj}(R)^\cong$ is the symmetric monoidal groupoid of projective, finitely generated $R$-modules for a discrete ring $R$, with symmetric monoidal structure given by the coproduct, the realization $| S^{-1} S |$ is equivalent to the $K$-theory space of $R$.
\item Suppose $M$ is a homotopy commutative topological monoid. Then for a given $m \in M$ one can consider the telescope
$$ M_m = \colim ( M \xrightarrow{m} M \xrightarrow{m} M \xrightarrow{m} \cdots ).$$
By iterating this construction over a generating set of $\pi_0 M$, one can construct a space $M_\infty$. Further, taking a plus construction one obtains a space $(M_\infty)^+$, together with the map
$$ M \rightarrow M_\infty \rightarrow (M_\infty)^+$$
which is a model for the group completion of $M$. This construction has been more or less considered by Quillen and has later been formalized by McDuff and Segal in \cite{McDuff1976} in the case when only a single generator $m$ is needed. A framework on how to deal with an arbitrary set of generators has been discussed by Randal-Williams \cite{10.1093/qmath/hat024} (See also \cite{Nikolaus_2017}).
\item A first definition of what would now be called an $\mathbb{E}_\infty$-monoid has been provided by Segal \cite{SEGAL1974293} with his notion of a $\Gamma$-space. Segal constructs a functor $B$ that deloops such a $\Gamma$-space. The group completion of a $\Gamma$-space $X$ can then be defined as $\Omega B X$.
\end{enumerate}

The purpose of this paper is twofold. The first is to put all these variations into one common framework, using the now established language of $\infty$-categories. Our model for what it means to be a ``sufficiently commutative higher algebraic gadget'' will be what is called an $\mathbb{E}_n$-monoid, also called $\mathbb{E}_n$-space, for $n \geq 2$. Roughly speaking, an $\mathbb{E}_n$-space is a coherently associative and unital monoid in spaces, that is furthermore commutative up to coherent homotopy, whereby $n$ determines how far this coherence holds. An $\mathbb{E}_n$-monoid is called grouplike if its set of path components, which naturally has the structure of a monoid, is a group. For purely formal reasons it holds that the inclusion of the full subcategory of the $\infty$-category of $\mathbb{E}_n$-monoids given by the grouplike objects has a left adjoint $(-)^{\gp}$, which is called \emph{group completion}. The key players will be three $\mathbb{E}_{n-1}$-monoidal $\infty$-categories one can construct from an $\mathbb{E}_n$-monoid $X$, called:
\begin{itemize}
\item $\B X$, the \emph{classifying} $\infty$-category of $X$,
\item $\E X$, the \emph{action} $\infty$-category of $X$, and
\item $\G X$, the \emph{double action} $\infty$-category of $X$.
\end{itemize}
The chosen terminology reflects that $\E X$ and $\G X$ are obtained as unstraightenings of the action of $X$ on itself, and on $X \times X$, respectively. These $\mathbb{E}_{n-1}$-monoidal $\infty$-categories sit together in a pullback square
$$\xymatrix{
\G X \ar[r]^{p_+} \ar[d]_{p_-} & \E X  \ar[d]^p \\
\E X \ar[r]_p & \B X.
}$$
with the pullback taken in the $\infty$-category of $\mathbb{E}_{n-1}$-monoidal $\infty$-categories.

The central theorem of this paper is the following. Let $\Spc$ be the $\infty$-category of spaces, and $\Catinf$ be the $\infty$-category of (small) $\infty$-categories.
\begin{theorem} \label{centraltheorem}
Let $n \geq 2$, or $n = \infty$. There exist functors
\emph{$$ \B, \E, \G : \EMon{n} \rightarrow \EMon{n-1}(\Catinf).$$}
such that for a given $\mathbb{E}_n$-monoid $X$, there are functors
\emph{$$\begin{array}{l}
X^\circlearrowleft : \B X \rightarrow \Spc, \\
X^\rightarrow : \E X \rightarrow \Spc
\end{array}$$}
and natural equivalences
\emph{$$ X^{\gp} \simeq \Omega | \B X | \simeq | \G X | \simeq \colim_{ \E X } X^\rightarrow \simeq \colim_{ \B X } X^\circlearrowleft \times X^\circlearrowleft $$}
of $\mathbb{E}_{n-1}$-monoids.\footnote{Here $|\mathcal{C}|$ refers to the realization of an $\infty$-category $\mathcal{C}$, i.e.\ the space obtained by inverting all arrows in $\mathcal{C}$.}
\end{theorem}

The connections to the previously mentioned constructions are as follows.
\begin{itemize}
\item A symmetric monoidal groupoid $S$ can be considered as an $\mathbb{E}_\infty$-monoid. Quillen's $S^{-1}S$ is then simply the homotopy category of the $\infty$-category $\G S$. If $S$ satisfies the additional assumption that for all $ n, s \in S$ the translation action  $s+: \text{Aut}(n) \rightarrow \text{Aut}(s+n)$  is injective, then $\G S$ is itself a $1$-category and thus equivalent to $S^{-1} S$. The generalization to arbitrary $\mathbb{E}_n$-monoids allows for a potential extension of the techniques Quillen introduced, for example in the context of topological $K$-theory or $K$ and $L$-theory of non-discrete ring spectra.
\item The $\infty$-category $\E X$ is a universal recipient for telescope constructions on $X$. We will show in Section \ref{telescopes} that there is a filtered diagram $D$ together with a functor $\iota : D \rightarrow \E X$ such that the comparison map
  $$ X_\infty = \colim_D (X^\rightarrow\circ \iota) \rightarrow \colim_{\E X} X^\rightarrow $$
models the plus construction, at least for the case when $X$ is an $\mathbb{E}_\infty$-monoid.
\end{itemize}




The second purpose of this paper will be a selection of applications of Theorem \ref{centraltheorem}. For example, we will give a construction of $\Omega X^{\gp}$ as a colimit.

\begin{theorem}[See Section \ref{loops}] \label{loopspaceascolimit}
Let $X$ be an $\mathbb{E}_n$-monoid for $n \geq 2$. There is a natural equivalence of $\mathbb{E}_1$-groups
$$ \Omega X^{\gp} \simeq \colim_{x \in \E X } \mathrm{Aut}(x) $$
where the colimit is taken in the $\infty$-category of $\mathbb{E}_1$-groups.
\end{theorem}

\begin{remark} Under the equivalence of $\mathbb{E}_1$-groups with pointed connected spaces, and the fact that weakly contractible colimits of pointed connected spaces are computed as colimits in (unpointed) spaces, we see that the above theorem is equivalent to
$$ X^{\gp}_1 \simeq \colim_{x \in \E X } B \mathrm{Aut}(x), $$
with the colimit taken in the $\infty$-category of spaces, and $X^{\gp}_1$ referring to the path component of the unit $1 \in X^{\gp}$.

If $E$ is a spectrum, we thus see that
$$ E[X^{\gp}_1] \simeq \colim_{x \in \E X } E[B \mathrm{Aut}(x)]$$
in the $\infty$-category of spectra, which recovers a version of the classical group completion theorem as a statement about homology. (See also Theorem \ref{homology}.) In this case, the colimit in question can be replaced by a filtered diagram, obtained as a telescope over chosen generators of $\pi_0 X$. This result is commonly used in work on homological stability, to obtain comparison results relating the homology of $X^{\gp}_1$ to the homology of the classifying spaces $B \mathrm{Aut}(x)$. See for example \cite[Section 3]{SzymikWahl2019} and \cite[Section 5]{Li2025}. In this sense, we would like to think of Theorem \ref{loopspaceascolimit} as a precursor in the land of unstable homotopy theory.
\end{remark}

Further, we will provide a relative version of the constructions provided above. Given a map $f : X \rightarrow Y$ of $\mathbb{E}_\infty$-monoids there are symmetric monoidal $\infty$-categories $\E (f)$ and $\G (f)$ such that the following theorem holds.

\begin{theorem}[See Section \ref{relative}]
If $f : X \rightarrow Y$ is a map of $\mathbb{E}_\infty$-monoids, there are natural cofiber sequences
$$ X \rightarrow Y \rightarrow |\E (f)|$$
and
$$ X^{\gp} \rightarrow Y^{\gp} \rightarrow |\G (f)|$$
of $\mathbb{E}_\infty$-monoids and $\mathbb{E}_\infty$-groups respectively.
\end{theorem}

Using this we prove a generalized version of the cofinality theorem.

\begin{theorem}[See Theorem \ref{cofinalitytheorem}]
Let $\iota : X \rightarrow Y$ be a cofinal inclusion of \emph{$\mathbb{E}_n$}-monoids for $n \geq 2$. Then \emph{$\cof(\iota^{\gp}) : X^{\gp} \rightarrow Y^{\gp}$} is a discrete commutative monoid. In particular,
$$ \pi_i X^{\gp} \rightarrow \pi_i Y^{\gp} $$
is injective in degree $0$ and an isomorphism for $i > 0$. Here, the cofiber $\cof(\iota^{\gp})$ is taken in the $\infty$-category of \emph{$\mathbb{E}_n$}-groups.
\end{theorem}

We note that a map of $\mathbb{E}_n$-monoids $\iota : X \rightarrow Y$ is in particular cofinal, if $\iota : X \rightarrow Y$ corresponds to an inclusion of a choice of path-components of $Y$, and the classical condition that for all $y \in Y$ there exists $y' \in Y$ and $x \in X$ such that $y \cdot y' \sim \iota(x)$ holds, compare \cite[p.\ 225]{10.1007/BFb0080003}.

\begin{remark}
It seems natural to ask if the main statements presented in this article can be generalized to the case where  $X$ is an internal $\mathbb{E}_n$-monoid in an arbitrary $\infty$-topos $\mathcal{X}$ rather than the $\infty$-topos of spaces. The main constructions of $\B X, \E X$ and $\G X$ (which should be understood as \emph{internal categories} in $\mathcal{X}$), as well as their basic properties certainly work in this more general situation. However, the current argument for the proof of Theorem \ref{Liscompletion}, which is the main motivation for considering $\G X$, relies on a combination of the Ore condition as well as a homology computation. The author does not know, at the time of writing, whether these arguments can be generalized.
\end{remark}

\subsection{Acknowledgements}

Very special thanks go to Maxime Ramzi, without whom this project would not have been possible. The author would also like to thank Christoph Winges, Thomas Nikolaus, Thomas Blom, Holger Reich, Lars Hesselholt, Vittorio di Fraia and David Kern for helpful discussions and commentary, as well as the patient anonymous referee, whose comments have significantly improved the quality of this article.

\section{Preliminaries}

We will use the language of $\infty$-categories and higher algebra as developed in \cite{luriehtt} and \cite{lurieha}, see also \cite{LandIntroductionInfinityCategories}. The $\infty$-category of spaces, also called animae or $\infty$-groupoids will be denoted by $\Spc$. The $\infty$-category of small $\infty$-categories will be denoted by $\Catinf$. The $\infty$-category of locally small $\infty$-categories will be denoted by $\widehat{\Catinf}$, where with \emph{locally small} we mean $\infty$-categories $\mathcal{C}$ with a potential class of objects, but such that the mapping spaces $\mathrm{Map}_\mathcal{C}(c,d)$ remain small for every pair of objects $c,d \in \mathcal{C}$. The $\infty$-category of presentable $\infty$-categories with maps given by left adjoint functors will be denoted by $\PrL$, and is a non-full subcategory of $\widehat{\Catinf}$. Spaces, viewed as $\infty$-groupoids, naturally embed into $\infty$-categories as a full subcategory. The inclusion $ \Spc \rightarrow \Catinf$ has a left adjoint $\real{-}$, called realization, and a right adjoint $\core$. The realization of an $\infty$-category $\mathcal{C}$ can be understood as the localization at all morphisms, i.e.\ $|\mathcal{C}| = \mathcal{C}[\mathcal{C}^{-1}]$. The $\infty$-groupoid $\core(\mathcal{C})$ can be defined as the maximal sub-$\infty$-groupoid of $\mathcal{C}$.  A crucial tool will be the straightening/unstraightening equivalence due to Lurie.

\begin{theorem}[\cite{luriehtt}, Chapter 3.2, see also \cite{LandIntroductionInfinityCategories}, Theorem 3.3.16.]
For every $\infty$-category $\mathcal{C}$ there is an equivalence
$$ \mathrm{Str} : \mathrm{LFib}( \mathcal{C} ) \simeq \Fun( \mathcal{C}, \Spc ) : \mathrm{Un} $$
where $\mathrm{LFib}( \mathcal{C} )$ denotes the full subcategory of the slice $\Catinf{/ \mathcal{C}}$ spanned by left fibrations.
\end{theorem}

The colimit of a functor $F : \mathcal{C} \rightarrow \Spc$ can be computed as a realization, as $| \text{Un}( F )| \simeq \colim_{\mathcal{C}} F$. Given $c \in  \mathcal{C}$, the mapping space functor $\text{Map}_\mathcal{C}(c,-) : \mathcal{C} \rightarrow \Spc$ corresponds to the left fibration $\mathcal{C}_{c /} \rightarrow \mathcal{C}$ of the undercategory of $c$. (See \cite{LandIntroductionInfinityCategories}, section 4.2.)

A functor $F : \mathcal{C} \rightarrow \mathcal{D}$ is called \emph{cofinal} if for all functors $G : \mathcal{D} \rightarrow \Spc$, the natural comparison map
$$ \colim_\mathcal{C} G \circ F \rightarrow \colim_\mathcal{D} G $$
is an equivalence. Note that by unstraightening, this in particular implies that $F$ induces an equivalence on the realizations $|\mathcal{C}| \simeq |\mathcal{D}|$. Quillen's Theorem A, see \cite{luriehtt}, proposition 4.1.3.1., states that $F$ is cofinal iff for all $d \in \mathcal{D}$, the $\infty$-category $d / F$, defined via the pullback
$$\xymatrix{
d / F \ar[r] \ar[d] & \mathcal{D}_{d /} \ar[d] \\
\mathcal{C} \ar[r]^F & \mathcal{D}.
}$$
of $\infty$-categories is weakly contractible. We would like to state the following alternate formulation of Quillen's Theorem A.

\begin{lemma} \label{quillensA}
A functor $F : \mathcal{C} \rightarrow \mathcal{D}$ is cofinal iff for all $d \in \mathcal{D}$,
$$ \colim_{c \in \mathcal{C}} \Map_\mathcal{D}( d, F(c) ) \simeq \pt. $$
\end{lemma}

\begin{proof}
The left fibration $\mathcal{D}_{d /}$ classifies the functor $\Map_\mathcal{D}( d, - )$. This implies that the functor $\Map_\mathcal{D}( d, F(-) )$ is classified by the left fibration $ d / F \rightarrow \mathcal{C}$ coming from the defining pullback
$$\xymatrix{
d / F \ar[r] \ar[d] & \mathcal{D}_{d /} \ar[d] \\
\mathcal{C} \ar[r] & \mathcal{D}.
}$$
But this implies that the colimit in question is given by the realization of $ d / F $, and hence the statement is equivalent to Quillen's Theorem A.
\end{proof}

We call a left adjoint functor $L : \mathcal{C} \rightarrow \mathcal{C}_0$ a \emph{Bousfield localization} if the right adjoint to $L$ is fully faithful. In this case $\mathcal{C}_0$ can be identified with the subcategory of $S$-local objects of $\mathcal{C}$, where $S$ is the class of morphisms of $\mathcal{C}$ that are sent to equivalences under $L$. Conversely, if $C$ is presentable and $S$ is a (small) set of morphisms, there exists a left adjoint $L$ to the inclusion of the subcategory of $C$ spanned by the $S$-local objects. For more information, see \cite{luriehtt}, section 5.2.7 and proposition 5.5.4.15.

\subsection{$\mathbb{E}_n$-monoids and $\mathbb{E}_n$-groups}

The theory of $\mathbb{E}_n$-monoids can be described with two equivalent approaches: algebraic theories and operads. Since the first are conceptually simpler to deal with we will use them here. We call an $\infty$-category $\mathcal{C}$ with finite products \emph{cartesian}. 

\begin{definition}[See also \cite{joyal2008}, \cite{cranch2010algebraic}, \cite{berman2019higher}, and \cite{Gepner_2015}, Appendix B]
An \emph{algebraic theory}, also called \emph{Lawvere theory}, is a cartesian $\infty$-category $\mathcal{L}$ together with a given object $x \in \mathcal{L}$ that generates $\mathcal{L}$ under finite products. If $\mathcal{C}$ is a cartesian category, a \emph{model} $M$ for $\mathcal{L}$ is a finite-product preserving functor $\mathcal{L} \rightarrow \mathcal{C}$. Write $\text{Mod}_\mathcal{L}(\mathcal{C})$ for the full subcategory of $\Fun(\mathcal{L}, \mathcal{C})$ consisting of the models.
\end{definition}

We define the following basic algebraic theories:
\begin{itemize}
\item We define $\mathcal{L}_{\mathbb{E}_0}$ as the opposite of the category $\Fin_*$ of finite pointed sets. A model for $\mathcal{L}_{\mathbb{E}_0}$ in a cartesian $\infty$-category $\mathcal{C}$ is the same as a pointed object.
\item We define $\mathcal{L}_{\mathbb{E}_1}$ as the opposite of the category of free, finitely generated monoids, viewed as a full subcategory of the $1$-category of monoids. Equivalently, $\mathcal{L}_{\mathbb{E}_1}$ has objects indexed by $\mathbb{N}$, and $\text{Hom}(n,1)$ given by the set of words in $n$ letters.\footnote{There is a canonical inclusion $\Delta^{\op} \rightarrow \mathcal{L}_{\mathbb{E}_1}$. The reader is invited to show that for a given cartesian $\infty$-category $\mathcal{C}$, restricting along this inclusion provides an equivalence between $\mathcal{L}_{\mathbb{E}_1}$-models and simplicial diagrams satisfying a Segal condition.}
\end{itemize}

Given two algebraic theories $\mathcal{L}$ and $\mathcal{L}'$, there is a Lawvere theory $\mathcal{L} \otimes \mathcal{L}'$, together with finite product preserving functors $\mathcal{L} \rightarrow \mathcal{L} \otimes \mathcal{L}'$ and $\mathcal{L}' \rightarrow \mathcal{L} \otimes \mathcal{L}'$, uniquely characterized by the property that for any cartesian $\infty$-category $\mathcal{C}$, there is a natural equivalence
$$ \text{Mod}_{\mathcal{L} \otimes \mathcal{L}'}( \mathcal{C} ) \simeq \text{Mod}_{\mathcal{L}} ( \text{Mod}_{\mathcal{L}'} ( \mathcal{C} ) ), $$
see \cite{berman2019higher}, section 3.
\begin{definition}
Define $\mathcal{L}_{\mathbb{E}_n}$ inductively as $\mathcal{L}_{\mathbb{E}_1} \otimes \mathcal{L}_{\mathbb{E}_{n-1}}$.  Define $\mathcal{L}_{\mathbb{E}_\infty} = \colim_{n} \mathcal{L}_{\mathbb{E}_n}$. Given a cartesian $\infty$-category $\mathcal{C}$, we write $\EMon{n}(\mathcal{C})$ for the $\infty$-category of models of $\mathcal{L}_{\mathbb{E}_n}$. If $\mathcal{C} = \Spc$ we simply write $\EMon{n}$ and call a model an $\mathbb{E}_n$-\emph{monoid} or $\mathbb{E}_n$-\emph{space}. In the case of $\mathcal{C} = \Catinf$ we call a model an $\mathbb{E}_n$-monoidal $\infty$-category and a map of models a strong $\mathbb{E}_n$-monoidal functor.
\end{definition}

\begin{remark}
We note that despite $\mathcal{L}_{\mathbb{E}_1}$ being a $1$-category, $\mathcal{L}_{\mathbb{E}_n}$ will not be a $1$-category any more for $n > 1$. For the case $n = \infty$, the algebraic theory $\mathcal{L}_{\mathbb{E}_\infty}$ is equivalent to the $(2,1)$-category $\text{Span}(\Fin)$, see \cite{cranch2010algebraic}.
\end{remark}

Alternatively, the category of $\mathbb{E}_n$-monoids can be defined using the $\mathbb{E}_n$-operad $\mathbb{E}_n^\otimes$, see \cite{lurieha}. This allows for a more flexible definition of models with respect to an $\mathbb{E}_n$-monoidal $\infty$-category $(\mathcal{C},\otimes)$, which does not need to be cartesian. We will call such models $\mathbb{E}_n$-algebras. The structure of an $\mathbb{E}_n$-monoidal $\infty$-category can be defined alternatively as a certain functor $\mathcal{C}^\otimes \rightarrow \mathbb{E}_n^\otimes$. An $\mathbb{E}_n$-lax monoidal functor between $\mathbb{E}_n$-monoidal $\infty$-categories will be a map of operads over $\mathbb{E}_n^\otimes$. Lax monoidal functors preserve algebras. We will later on need the statement that if a strong $\mathbb{E}_n$-monoidal functor $\mathcal{C} \rightarrow \mathcal{D}$ has a right adjoint $R$, the right adjoint will be lax monoidal.

Now let $X$ be an $\mathbb{E}_n$-monoid for $n \geq 1$, or $n= \infty$. Let $\B X$ be the one-object $\infty$-category with base-point $*$ and mapping space $X$. We call $\B X$ the \emph{classifying $\infty$-category} of $X$. This construction defines a functor, whose main properties are captured by the following proposition.

\begin{proposition}[See \cite{Gepner_2015}, Corollary 6.3.11.] \label{classifyinginftycat}
There is a functor
$$\B : \EMon{n} \rightarrow \EMon{n-1}(\Catinf ),$$
characterized by the property that $\B X$ has a single object and mapping space given by $X$. It is fully faithful and preserves products, and admits a right adjoint 
$$\mathbf{ \Omega } : \EMon{n-1}(\Catinf ) \rightarrow \EMon{n},$$
which maps an $\mathbb{E}_{n-1}$-monoidal $\infty$-category $\mathcal{C}$ to the $\mathbb{E}_n$-monoid $\mathbf{ \Omega } \mathcal{C}$ with underlying space $\mathrm{End}(1_\mathcal{C})$. The functor $\mathbf{ \Omega }$ is lax monoidal. The essential image of $\B$ is given by those $\mathbb{E}_{n-1}$-monoidal $\infty$-categories $\mathcal{C}$, such that the inclusion of the basepoint $\pt \rightarrow \mathcal{C}$ is essentially surjective.
\end{proposition}

\begin{remark}
Let us use the notation $BX$ for the \emph{classifying space} of an $\mathbb{E}_n$-monoid $X$, which is by definition the realization $| \B X |$. We would like to stress the distinction between $\B X $ and $BX$. In general, $\B X$ is not an $\infty$-groupoid, and hence not equivalent to its realization - this is the case iff $X$ is an $\mathbb{E}_n$-group.
\end{remark}

\begin{remark}
Let us point out that the adjunction $\B \dashv \Omega$ works in the greater generality where $\Catinf$ is replaced by the $\infty$-category $\widehat{\Catinf}$ of locally small $\infty$-categories, in the sense that there is an adjunction
\[\begin{tikzcd}
\mathbb{E}_n\mathrm{Mon}
  \arrow[r, bend left=20, "\mathbf{B}"]
&
\mathbb{E}_{n-1}\mathrm{Mon}(\widehat{\mathrm{Cat}_\infty})
  \arrow[l, bend left=20, "\Omega"]
\end{tikzcd}
\]
simply because $\Omega( \mathcal{C}, c) = \mathrm{Map}_\mathcal{C}(c,c)$ is a small space for locally small $\mathcal{C}$, and thus $\Omega$ is well-defined as a functor.
\end{remark}

Assume $n \geq 1$. Given an $\mathbb{E}_n$-monoid $X$, since the functor $\pi_0$ preserves finite products, the set $\pi_0 X$ naturally has the structure of a monoid (commutative if $n \geq 2$). We call $X$ \emph{grouplike} if $\pi_0 X$ is a group. Write $\EGrp{n}$ for the full subcategory of $\EMon{n}$ spanned by grouplike $\mathbb{E}_n$-monoids. Recall that there are two adjunctions
\[
\begin{tikzcd}
\Catinf \ar[rr,bend left,"|-|",""{name=A}] 
\ar[rr,bend right,"\rotadj","\core"',{name=C,below}]& & 
\Spc. \ar[ll,"\rotadj"{name=B,above}]
\end{tikzcd}
\]

\begin{proposition}
The inclusion $\EGrp{n} \hookrightarrow \EMon{n}$ has left and right adjoints, called \emph{group completion} $(-)^{\gp}$ and \emph{units} $(-)^\times$, respectively. Concretely, they are constructed as
$$ (X)^{\gp} = \Omega | \B X | $$
and
$$ X^\times = \Omega \core( \B X ).$$
Moreover, $(-)^{\gp}$ preserves finite products.
\end{proposition}

\begin{remark}
Using the notation $B X$ for $| \B X |$, we obtain the familiar formula $ X^{\gp} \simeq \Omega B X$.
\end{remark}

\begin{proof} Let $h : \Catinf \rightarrow \Cat$ be the functor that assigns to an $\infty$-category its homotopy category. There is the commutative square
$$\xymatrix{
\EMon{n} \ar[r]^{\B} \ar[d]^{\pi_0} & \EMon{n-1}(\Catinf) \ar[d]^h \\
\EMon{n}(\Set)       \ar[r]^{\B} & \EMon{n-1}(\Cat),
}$$
from which we see that $X$ is grouplike iff $h \B X$ is a groupoid iff $\B X$ is an $\infty$-groupoid.

The two pairs of adjoint functors
\[
\begin{tikzcd}
\Catinf \ar[rr,bend left,"|-|",""{name=A}] 
\ar[rr,bend right,"\rotadj", "\core"',{name=C,below}]& & 
\Spc \ar[ll,"\rotadj"{name=B,above}]
\end{tikzcd}
\]
give the adjunctions
\[
\begin{tikzcd}
\EMon{n-1}(\Catinf) \ar[rr,bend left,"|-|",""{name=A}] 
\ar[rr,bend right,"\rotadj", "\core"',{name=C,below}]& & 
\EMon{n-1}(\Spc) \ar[ll,"\rotadj"{name=B,above}]
\end{tikzcd}
\]
since both $|-|$ and $\core$ preserve finite products. Both functors furthermore preserve the condition for an $\mathbb{E}_{n-1}$-monoidal $\infty$-category $\mathcal{C}$ that $ \pt \rightarrow \mathcal{C}$ is essentially surjective, and hence by Proposition \ref{classifyinginftycat} induce the pair of adjoint functors
\[
\begin{tikzcd}
\EMon{n} \ar[rr,bend left,"(-)^{\gp}",""{name=A}] 
\ar[rr,bend right,"\rotadj","(-)^\times"',{name=C,below}]& & 
\EGrp{n}. \ar[ll,"\rotadj"{name=B,above}]
\end{tikzcd}
\]
The statement that $(-)^{\gp}$ preserves finite products follows from the analogous statement about $|-|$.
\end{proof}

For the following corollary, write $(-)^{\gp,n} : \EMon{n} \rightarrow \EGrp{n}$ for the group completion functor on $\mathbb{E}_n$-monoids.

\begin{corollary} \label{groupcompletionforget}
Let $U_n : \EMon{n} \rightarrow \EMon{n-1}$ be the forgetful functor induced by the functor $\mathcal{L}_{\EMon{n-1}} \rightarrow \mathcal{L}_{\EMon{n}}$. Then there is a commutative square
$$\xymatrix{
\EMon{n} \ar[r]^{U_n} \ar[d]^{(-)^{\gp},n} & \EMon{n-1} \ar[d]^{(-)^{\gp},n-1} \\
\EGrp{n} \ar[r]^{U_n} & \EGrp{n-1}.
}$$
\end{corollary}

\begin{proof}
Since
$$(-)^{\gp,n-1} : \EMon{n-1} \rightarrow \EGrp{n-1} $$
preserves finite products, it induces the functor
$$ \EMon{n} \rightarrow \EMon{1}( \EGrp{n-1} ) \cong  \EGrp{n},$$
which can be seen to be equivalent to the group completion $(-)^{\gp,n}$, since it is by construction left adjoint to the inclusion of $\EGrp{n}$ into $\EMon{n}$. This in particular implies commutation with the forgetful functors.
\end{proof}

\begin{remark}
If $\mathcal{C}$ is an $\infty$-groupoid, the condition that $ \pt \rightarrow \mathcal{C}$ is essentially surjective becomes equivalent to $\pi_0 \mathcal{C} = 0$. It is now straightforward to see that the functor $\B$ induces an equivalence
$$\EGrp{1} \xrightarrow{\simeq} \Spc_{\geq 1}$$
where the righthand side is the $\infty$-category of connected spaces. Further iterating this gives the equivalence
$$\EGrp{n} \xrightarrow{\simeq} \Spc_{\geq n}$$
for any $n$, thus recovering the celebrated May recognition principle.
\end{remark}

\begin{remark}
It can be shown that the category of $\mathbb{E}_1$-groups is again the category of models for an algebraic theory, namely $\mathcal{L}_{\EGrp{1}} = \mathcal{L}_{\text{Grp}}$ is the category representing the $1$-categorical theory of groups, i.e.\ it is equivalent to the opposite of the category of free groups. The free group on a single generator, $\mathbb{Z}$, has the inverse map $(-1) : \mathbb{Z} \rightarrow \mathbb{Z}$. Let $X$ be an $\mathbb{E}_1$-group. Viewed as a functor $\mathcal{L}_{\text{Grp}} \rightarrow \Spc$ it is clear that there exists a homotopy inverse $\iota$ to the multiplication, given as the image of $(-1)$. It follows that the path components of the underlying space of $X$ must all be equivalent, hence
$$X \simeq \pi_0(X) \times X_1,$$
in the $\infty$-category of spaces, where $X_1$ is the path component of the identity.
\end{remark}

\subsection{The Plus construction as a localization}

We would like to collect a few notions around the plus construction. Recall that a group $H$ is called \emph{perfect}, if $H_{ab} = 0$. A group $G$ is called \emph{hypoabelian} if there does not exist a non-trivial perfect subgroup. We call a space $X$ hypoabelian if $\pi_1(X,x)$ is \emph{hypoabelian} for all choices of $x$. Denote by $\Spc^{\text{hypo}}$ the full subcategory of $\Spc$ spanned by hypoabelian spaces.

\begin{theorem} \label{plusconstruction}
The inclusion of \emph{$\Spc^{\text{hypo}}$} into \emph{$\Spc$} has a left adjoint, given by the plus construction $(-)^+$. It is the Bousfield localization with respect to the class of acyclic maps.
\end{theorem}

This theorem is discussed in a classical form in \cite{BERRICK1999467}. A modern proof can be found in \cite{HoyoisPlusConstruction}.

Note that if $X$ is an $\mathbb{E}_1$-group, its underlying space is hypoabelian, and thus local for the class of acyclic maps. This is clear for path connected $\mathbb{E}_1$-groups $X_1$, as the fundamental group of such an $X_1$ is abelian, and it follows for general $\mathbb{E}_1$-groups since all path components are mutually equivalent. Note that the grouplike assumption is necessary. The $\mathbb{E}_1$-monoid $\Fin^\cong$ has the property that $\pi_1(\Fin^\cong, \emptyset) = 0$ is abelian, but not in general for the rest of the path components, which are given by $B\Sigma_n$'s, which are not hypoabelian for $n > 4$. The following is immediate from Theorem \ref{plusconstruction}.

\begin{corollary} \label{pluscorollary}
Let $X \rightarrow Y$ be an acyclic map of spaces, and suppose $Y$ can be equipped with an $\mathbb{E}_1$-group structure. Then $Y$ is a model for the plus construction of $X$. Moreover, a map $f : X \rightarrow Y$ between spaces $X$ and $Y$ which can both be equipped with $\mathbb{E}_1$-group structures is acyclic iff it is an equivalence.
\end{corollary}

\section{The action $\infty$-category of an $\mathbb{E}_n$-monoid}

Assume that $G$ is a discrete group. There are two groupoids one can associate to $G$. One is the groupoid $BG$ which has a single object $*$ and $\text{hom}(*,*) = G$. The other is called $EG$, which has as objects the underlying set $G$, and morphisms $ g \xrightarrow{h} g'$ whenever $h g = g'$. Since such an $h$ is necessarily unique, every object in $EG$ is both initial and terminal, and hence the realization of $EG$ is contractible. There is a canonical functor
$ EG \rightarrow BG $ which sends an arrow $ g \xrightarrow{h} g'$ to $h \in G$. It is a standard exercise in topology that the realization of this functor is a map of topological spaces is a model for the universal cover of the classifying space of $G$. We want to generalize this construction.

Let $X$ be an $\mathbb{E}_n$-monoid for $n \geq 1$, or $n= \infty$ and let $\B X$ be the one-object $\infty$-category with base-point $*$ and mapping space $X$, as discussed in the previous section.\footnote{We should remark that in the following, starting from Definition \ref{actioncategorydefinition} until before Proposition \ref{monoidalstructure}, only the existence of an $\mathbb{E}_1$-structure will be used, and the additional $\mathbb{E}_n$-structure will be temporarily forgotten.} There is a functor
$$ X^\circlearrowleft : \B X \rightarrow \Spc$$
given by the action of $X$ on itself via left multiplication, or equivalently said, the functor corepresented by the base point.
\begin{definition} \label{actioncategorydefinition}
We define $p : \E X \rightarrow \B X$ as the left fibration classifying $X^\circlearrowleft$, or equivalently as the projection from the under-category $\B X_{*/} \rightarrow \B X$. We call $\E X$ the \emph{action $\infty$-category} of $X$.
We denote by $X^\rightarrow$ the functor given by the composition
$$X^\rightarrow : \E X \xrightarrow{p} \B X \xrightarrow{ X^\circlearrowleft } \Spc. $$
We define $\G X$ as the pullback in $\Catinf$,
$$\xymatrix{
\G X \ar[r]^{p_+} \ar[d]_{p_-} & \E X  \ar[d]^p \\
\E X \ar[r]_p & \B X.
}$$
We call $\G X$ the \emph{double action $\infty$-category} of $X$. 
\end{definition}

We note that by construction, $\G X \rightarrow \E X$ classifies the functor $X^\rightarrow : \E X \rightarrow \Spc$, and $\G X \rightarrow \B X$ classifies the functor $X^\circlearrowleft \times X^\circlearrowleft: \B X \rightarrow \Spc$. Let us give a concrete description of the $\infty$-categories $\E X$ and $\G X$. The fiber of $\E X \rightarrow \B X$ over the single object of $\B X$ is given by $X$, hence we see that the objects of $\E X$ are given by the objects of the $\infty$-groupoid $X$. Similarly, the objects of $\G X$ are given by the objects of $X \times X$. For a given $x \in X$ we also have the map given by right multiplication $\cdot x : X \rightarrow X$. We claim that the mapping spaces of $\E X$ are given as the fibers of this map.

\begin{lemma} \label{mappingspaces}
For $x,y \in \E X$ there is a natural pullback square of spaces
$$\xymatrix{
\Map_{\E X}( x,y) \ar[r] \ar[d] & X \ar[d]^{\cdot x} \\
\pt \ar[r]^{y} & X.
}$$
For $(x_1,x_2), (y_1,y_2) \in \G X$ there is a natural pullback square of spaces
$$\xymatrix{
\Map_{\G X}(( x_1,x_2), (y_1,y_2)) \ar[r] \ar[d] & X \ar[d]^{(\cdot x_1, \cdot x_2)} \\
\pt \ar[r]^{(y_1, y_2)} & X \times X.
}$$
\end{lemma}

\begin{proof}
Since $\E X = \B X_{*/}$ is an under-category, by the usual formula for mapping spaces in under-categories we have the pullback
$$\xymatrix{
\Map_{\E X}( x,y) \ar[r] \ar[d] & \Map_{\B X}( *, *) \ar[d]^{x^*} \\
\pt \ar[r]^y & \Map_{\B X}( *, *).
}$$
Under the canonical identification of $\Map_{\B X}( *, *) = X$, the statement follows. For the mapping spaces of $\G X$, since $\G X$ is defined as a pullback of $\infty$-categories we have the pullback square
$$\xymatrix{
\Map_{\G X}(( x_1,x_2), (y_1,y_2)) \ar[r] \ar[d] & \Map_{\E X}( x_1, y_1) \ar[d] \\
\Map_{\E X}( x_2, y_2) \ar[r] & \Map_{\B X}( *, *).
}$$
One can paste this pullback with the previous square twice and again identify $\Map_{\B X}( *, *)$ with $X$ to get the diagram
$$\xymatrix{
\Map_{\G X}(( x_1,x_2), (y_1,y_2)) \ar[r] \ar[d] & \Map_{\E X}( x_1, y_1) \ar[d] \ar[r] & \pt \ar[d]^{y_1}\\
\Map_{\E X}( x_2, y_2) \ar[r] \ar[d] & X \ar[r]^{\cdot x_1} \ar[d]^{\cdot x_2} & X  \\
\pt \ar[r]^{y_2} & X & 
}$$
where all squares are pullbacks in $\Spc$. Thus $\Map_{\G X}(( x_1,x_2), (y_1,y_2))$ is the limit over the bottom right boundary of the total diagram. The claimed pullback square for $\Map_{\G X}(( x_1,x_2), (y_1,y_2))$ immediately follows by exchanging the order in which to do the limit.
\end{proof}

Before we continue, a few words on the action category $\E X$. As observed before, objects $x \in \E X$ can be identified with the objects of $X$, and by the previous lemma, morphisms $f : x \rightarrow y$ are given by $t \in X$ together with a path $tx \sim y$. The functor $X^\rightarrow : \E X \rightarrow \Spc$ is given on objects and morphisms by
$$\begin{tikzcd}[row sep=-1pt]
x  &[-20pt] \mapsto &[-20pt] &[-20pt] X &[-20pt] \\
(f : x \xrightarrow{t} y) & \mapsto & X \ar[rr, "t \cdot"]  & & X.
\end{tikzcd}$$
There is a natural morphism $ 1 \rightarrow x $ for any $x \in X$, given by $x$ acting on $1$. It will come to no surprise that this makes $1$ into an initial object.

\begin{corollary} \label{initialobject}
The object $1 \in \E X$ is initial.
\end{corollary}

\begin{proof}
This is immediate, since $\cdot 1 \simeq \text{id}_X$, hence $\Map_{\E X}( 1,y) \simeq \pt$ for all $y \in X$. 
\end{proof}

\begin{example} If $X$ is a discrete monoid, the $\infty$-category $\E X$ is a $1$-category that can be concretely described as having as objects the set $X$, and a morphism $ x \xrightarrow{t} y$ whenever $tx = y$. The $\infty$-category $\G X$ is a $1$-category as well, with objects the set $X \times X$ and a morphism $(x_1, x_2) \xrightarrow{t} (y_1, y_2)$ whenever $tx_1 = y_1$ and $tx_2 = y_2$. \begin{enumerate}
\item As a special case, if $X = (\mathbb{N},+)$ then $\E\mathbb{N} = (\mathbb{N},\leq)$, viewed as a category. The category $\G\mathbb{N}$ is also given as a poset, with $(n_1, n_2) \leq (m_1, m_2)$ whenever $m_1 - n_1 = m_2 - n_2 \geq 0$. It should be immediate that the realization of $\G\mathbb{N}$ is $\mathbb{Z}$, viewed as a discrete space.
\item If $X = G$ is a discrete group, as  mentioned before, $\B G$ is a model for the classifying space of $G$, and $\E G$ is its universal cover (both represented by groupoids). The space $\G G$ may be less familiar. The reader is invited to show that it is a groupoid with $G$ many contractible components, each isomorphic to $\E G$.
\end{enumerate}
\end{example}

\begin{example} \label{finset}
The category of finite sets is a symmetric monoidal with tensor product given by the disjoint union. Taking the core of this category we get the $\mathbb{E}_\infty$-monoid of finite sets with bijections $(\Fin^\cong, \oplus )$. This symmetric monoidal category is the free $\mathbb{E}_\infty$-monoid on a single generator, \cite[p.\ 65]{cranch2010algebraic}. A calculation using Lemma \ref{mappingspaces} shows that $\E(\Fin^\cong)$ is a $1$-category and can be identified with the category $\text{FinInj}$ of finite sets with injections. The $\infty$-category $\G(\Fin^\cong)$ is a 1-category as well. Given the Barratt-Priddy-Quillen theorem, Theorem \ref{Liscompletion} implies that
$$ | \G(\Fin^\cong) | \simeq \Omega^\infty \mathbb{S}.$$
\end{example}

\begin{example}
More generally, suppose $(S,+)$ is a symmetric monoidal groupoid. We can view $S$ as a 1-truncated $\mathbb{E}_\infty$-monoid. Then $\B S$ will be a symmetric monoidal $(2,1)$-category. Assume that $S$ satisfies the following condition:
  \begin{equation}\label{eqn:condition}
    \forall n, s \in S: \text{The translation action } s+: \text{Aut}(n) \rightarrow \text{Aut}(s+n) \text{ is injective.}\tag{*}
  \end{equation}
This is the case for example for $S = \Fin^\cong$ or when $S = \text{Proj}(R)$, the groupoid of finitely generated projective $R$-modules for a discrete ring $R$. Then by Lemma \ref{mappingspaces} it is immediate from the long exact sequence in the homotopy groups of a fiber sequence that the $\infty$-categories $\E S$ and $\G S$ are $1$-categories. In fact, (\ref{eqn:condition}) holds iff $\E S$ is a 1-category. The categories $\E S$ and $\G S$ have already been considered by Quillen in \cite{10.1007/BFb0080003}. In \cite{10.1007/BFb0080003}, the category $\E S$ is called $\left\langle S, S \right\rangle$, and given concretely by
$$\left\langle S, S \right\rangle = \left\{ \begin{array}{ccc}
\text{Objects} &=& \text{Objects of } S \\
\text{Morphisms } x \rightarrow y &=& \text{Equivalence classes of } (k, \alpha)
\end{array}\right.$$
where $x,y,k \in S$ and $\alpha : k+x \cong y$ is a morphism in $S$. The equivalence relation on the morphisms corresponds exactly to the process of taking the homotopy category of $\E S$. More generally, without condition (\ref{eqn:condition}), the category $\left\langle S, S \right\rangle$ is just the homotopy category of $\E S$.

The category $\G S$ is referred to as $S^{-1} S$. It is given concretely by
$$S^{-1} S = \left\{\begin{array}{ccc}
 \text{Objects} &=& \text{Pairs of objects of } S \\
\text{Morphisms } (x_1, x_2) \rightarrow (y_1, y_2) &=& \text{Equivalence classes of } (k, \alpha, \beta)
\end{array}\right.$$
where $x_1, x_2,y_1, y_2,k \in S$ and $\alpha : k+x_1 \cong y_1, \beta : k+ x_2 \cong y_2$ are morphisms in $S$. Again, in the absence of condition (\ref{eqn:condition}), $S^{-1} S$ is equivalent to the homotopy category of $\G S$ and $(*)$ implies that $\G S$ is itself a $1$-category. Quillen proves that in the case of $S$ being $\text{Proj}(R)$ for a discrete ring, the realization $|S^{-1} S|$ is equal to the $K$-theory space of $R$. We will recover this result later on.
\end{example}

Since, by definition, the map $p_- : \G X \rightarrow \E X$ arises via pullback from the map $\E X \rightarrow \B X$, which classifies the action of $X$ on itself, it follows that the fiber of $p_-$ over the point $1 \in \E X$ agrees with the fiber of $\E X \rightarrow \B X$ over the canonical basepoint, in other words the underlying space of $X$. Let us refer to the resulting natural functor $\mathrm{pos}: X \rightarrow \G X$ obtained by the pullback
$$\xymatrix{
X \ar[r]^{\mathrm{pos}} \ar[d] & \G X \ar[d]^{p_-} \\
\pt \ar[r]^{1} & \E X,
}$$
taken in $\Catinf$, as the \emph{inclusion of positive elements}.

\begin{lemma} \label{grouplikesituation}
Suppose that $X$ is grouplike. Then $\B X$, $\E X$ and $\G X$ are $\infty$-groupoids, $\E X$ is contractible, and the map $\mathrm{pos} : X \rightarrow \G X$ is an equivalence of $\infty$-groupoids.
\end{lemma}

\begin{proof}
If $X$ is grouplike, then the homotopy category $h\B X \simeq \B\pi_0X$ is a groupoid, hence $\B X$ is an $\infty$-groupoid. Since $\E X \rightarrow \B X$ is a left fibration, this implies that $\E X$ is an $\infty$-groupoid as well. The same argument works for $\G X$. Due to Corollary \ref{initialobject}, the inclusion of the initial object $1 : \pt \rightarrow \E X$ is an equivalence, which also implies that $\mathrm{pos} : X \rightarrow \G X$ is an equivalence, as it is obtained via pullback from the inclusion of $1$ in $\E X$.
\end{proof}

\begin{lemma} \label{threeequivalentcolimits}
There are canonical equivalences
$$ | \G X | \simeq \colim_{\E X} ( X^\rightarrow ) \simeq  \colim_{\B X} ( X^\circlearrowleft \times X^\circlearrowleft ). $$
\end{lemma}

\begin{proof}
This is immediate by noting that the left fibration $\G X \rightarrow \E X$ classifies $X^\rightarrow$, and the left fibration $\G X \rightarrow \E X \rightarrow \B X$ classifies $X^\circlearrowleft \times X^\circlearrowleft.$
\end{proof}

We denote by $L(X)$ any of the three equivalent functors in \ref{threeequivalentcolimits}. We obtain a functor
$$ L : \EMon{n} \rightarrow \Spc. $$
Moreover, define a natural transformation $\eta_X : X \rightarrow L(X)$ as the composite
$$ X \xrightarrow{\mathrm{pos}} \G X \rightarrow | \G X |. $$
The following lemma is immediate from \ref{grouplikesituation}.
\begin{lemma} \label{etaequivalenceforgroups}
If $X$ is grouplike, then $\eta_X : X \rightarrow L(X)$ is an equivalence.
\end{lemma}

The $\infty$-categories $\B X$, $\E X$ and $\G X$ carry extra structure which will be important.

\begin{proposition} \label{monoidalstructure}
Let $X$ be an $\mathbb{E}_n$-monoid.
\begin{enumerate}
\item The $\infty$-categories $\B X$, $\E X$ and $\G X$ admit canonical
$\mathbb{E}_{n-1}$-monoidal structures.
\item The natural functors
\[
\E X \rightarrow \B X,
\qquad
p_{\pm} : \G X \rightarrow \E X
\]
admit canonical lifts to $\mathbb{E}_{n-1}$-monoidal functors.
\item The $\infty$-groupoid $L(X)=|\G X|$ admits a canonical
$\mathbb{E}_{n-1}$-monoid structure. Moreover, the functor $L$ refines to a functor
\[ L : \EMon{n} \rightarrow \EMon{n-1}. \]
\item The natural maps
\[
\mathrm{pos} : X \rightarrow \G X,
\qquad
\eta_X : X \rightarrow L(X)
\]
admit canonical lifts to $\mathbb{E}_{n-1}$-monoidal maps.
\end{enumerate}
\end{proposition}

We will henceforth, for a given $\mathbb{E}_n$-monoid $X$, implicitly equip $\B X$, $\E X$, $\G X$ as well as $L(X)$ with their canonical $\mathbb{E}_{n-1}$-monoidal structures.

\begin{proof}
If $X$ is an $\mathbb{E}_n$-monoid, Proposition \ref{classifyinginftycat} states that $\B X$ has a canonical $\mathbb{E}_{n-1}$-monoidal structure. It follows formally that the under category $\E X = \B X_{*/}$ is equipped with a canonical $\mathbb{E}_{n-1}$-monoidal structure such that the projection $\E X \rightarrow \B X$ also has the structure of a (strongly) $\mathbb{E}_{n-1}$-monoidal functor, see \cite[Theorem 2.2.2.4 and Remark 2.2.2.5]{lurieha}. Note that $\G X$ is obtained as the pullback of  $\E X \rightarrow \B X$ along itself in $\Catinf$. The forgetful functor $\EMon{n-1}( \Catinf ) \rightarrow \Catinf$ preserves limits \cite[Corollary 3.2.2.4]{lurieha}, so simply taking the pullback instead in $\EMon{n-1}( \Catinf )$ of the canonical $\mathbb{E}_{n-1}$-monoidal lift of $\E X \rightarrow \B X$ equips $\G X$ with a canonical $\mathbb{E}_{n-1}$-monoidal structure and gives the functors $p_{\pm} : \G X \rightarrow \E X$ canonical $\mathbb{E}_{n-1}$-monoidal lifts. This settles points (1) and (2).

For (3), the statement that $L(X) = | \G X |$ is canonically an $\mathbb{E}_{n-1}$-monoid follows from noting that realization $|-| : \Catinf \rightarrow \Spc$ commutes with finite products and is thus strong monoidal, hence induces a functor $$|-| : \EMon{n-1}( \Catinf ) \rightarrow \EMon{n-1}.$$ Point $(4)$ is clear for similar reasons, again since the functors in question are defined via pullback, so taking the defining pullbacks instead in $\EMon{n-1}( \Catinf )$ produces the lifts.
\end{proof}

\subsection{Actions and left modules}

If $f : X \rightarrow Y$ is a map of $\mathbb{E}_1$-monoids, we can think of $f$ as equipping $Y$ with an $X$-action via left multiplication. To make this observation precise, let us discuss actions of an $\mathbb{E}_1$-monoid in this section.

\begin{definition}
Let $X$ be an $\mathbb{E}_1$-monoid. A \emph{left $X$-action} on a space $Z$ is a functor $\underline{Z} : \B X \rightarrow \Spc$ such that $\underline{Z}(*) \simeq Z$. Define the $\infty$-category of left actions as $\Fun( \B X, \Spc )$.

Similarly, \emph{right $X$-actions} are defined as functors $\B X^{\op} \rightarrow \Spc$, and the $\infty$-category of right actions is defined as $\Fun( \B X^{\op}, \Spc )$.
\end{definition}

If $f : X \rightarrow Y$ as a map of $\mathbb{E}_1$-monoids, there is an induced adjunction on the $\infty$-categories of left actions
\[
\begin{tikzcd}
  \Fun( \B X, \Spc ) \arrow[r, bend left=10, "f_!"] &
  \Fun( \B Y , \Spc ) \arrow[l, bend left=10, "f^*", "\rotadj"']
\end{tikzcd}
\]
where $f^*$ is given by precomposition with $\B (f)$ and $f_!$ is given by left Kan extension.

The notion of left action for an $E_1$-monoid has a generalization to what is called a \emph{left module} for an algebra object. The theory of left modules for an algebra is developed in \cite[Chapter 4]{lurieha}. To make the relation to Lurie's more general framework understandable, let us sketch how the relation works.

For a given monoidal $\infty$-category $\mathcal{C}$, and $A$ an $\mathbb{E}_1$-algebra object, there is an associated $\infty$-category $$ \mathrm{LMod}_A(\mathcal{C})$$
of \emph{left modules} over $A$. If $\mathcal{C}$ is presentably monoidal, i.e.\ $\mathcal{C}$ is presentable and the tensor product on $\mathcal{C}$ commutes with colimits in both variables, then $ \mathrm{LMod}_A(\mathcal{C})$ is presentable as well, \cite[Corollary 4.2.3.7]{lurieha}. Moreover, it contains a natural choice of object, with $A \in \mathrm{LMod}_A(\mathcal{C})$ considered as a left module acting on itself. Structurally, the construction of the $\infty$-category of left $A$-modules can be understood in the following way.

\begin{theorem}[\cite{lurieha}, p.\ 750] Let $\mathcal{C}$ be a presentably monoidal $\infty$-category, and denote by $\mathrm{Alg}(\mathcal{C})$ the $\infty$-category of $\mathbb{E}_1$-algebras in $\mathcal{C}$. There exists a functor
\[\begin{array}{rcl}
\theta : \mathrm{Alg}(\mathcal{C}) & \rightarrow & \mathrm{LMod}_\mathcal{C}(\PrL)_{\mathcal{C} /} \\
                  A  & \mapsto & \mathrm{LMod}_A(\mathcal{C})
\end{array}  
\]
Moreover, $\theta$ is fully faithful and has a right adjoint $\mathrm{End}$.
\end{theorem}

The right adjoint $\mathrm{End}$ works the following way. If $i_! : \mathcal{C} \rightarrow \mathcal{M}$ is a $\mathcal{C}$-linear, left adjoint functor, then the image of the monoidal unit $1 \in \mathcal{C}$ gives an object $M = i_!(1)$ in $\mathcal{M}$. Then $\mathrm{End}( i_! )$ is defined to be the endomorphism algebra of $M$.

In the concrete case of $\mathcal{C} = \Spc$, the situation simplifies considerably. The presentable $\infty$-category $\Spc$ is the unit for the Lurie tensor product on $\PrL$, hence $\mathrm{LMod}_\Spc(\PrL) \simeq \PrL$. Algebras in $\Spc$ are simply $\mathbb{E}_1$-monoids, \cite[Chapter 4.1.2]{lurieha}, and so we obtain a functor
\[ \begin{array}{rcl}
 \theta : \EMon{1} & \rightarrow & \PrL_{ \Spc / } \simeq \EMon{0}(\PrL) \\
             X & \mapsto & \mathrm{LMod}_X(\Spc).
\end{array} \]
The $\infty$-category $\PrL_{ \Spc / }$ is equivalently described as the $\infty$-category of presentable $\infty$-categories together with a choice of object, and colimit preserving functors preserving the chosen objects. The right adjoint to $\theta$ sends a pair $(\mathcal{D}, d \in \mathcal{D})$ to
$$ \mathrm{End}_\mathcal{D}(d) = \mathrm{Map}_\mathcal{D}(d,d) $$
viewed as an $\mathbb{E}_1$-monoid, and thus identifies with the functor $\Omega$ given in Proposition \ref{classifyinginftycat}.

We can construct another functor $\EMon{1} \rightarrow \EMon{0}(\PrL)$, defined as the composite functor
\[
\EMon{1} \xrightarrow{\B} \EMon{0}(\Catinf)
\xrightarrow{P} \EMon{0}(\PrL).
\]
where $P(\mathcal{C}) = \Fun( \mathcal{C}^{\op}, \Spc )$. This composite takes an $\mathbb{E}_1$-monoid $X$ to the $\infty$-category of right actions $\Fun( \B X^{\op}, \Spc )$.

\begin{proposition}
There exists a natural equivalence
$$ P \circ \mathbf{B} \simeq \theta $$
of functors $\EMon{1} \rightarrow \EMon{0}(\PrL)$. In particular, for any given 
$X \in \EMon{1}$, there is a canonical equivalence
$$\Fun(\B X^{\op}, \Spc) \simeq  \mathrm{RMod}_X(\Spc).$$
\end{proposition}

\begin{remark}
Of course, the equivalence 
$$\Fun(\B X, \Spc) \simeq  \mathrm{LMod}_X(\Spc)$$
follows dually.
\end{remark}

\begin{proof}
Since left adjoints to a given functor are unique up to natural equivalence, it suffices to show that the functor $P \circ \B : \EMon{1} \rightarrow \EMon{0}(\PrL)$ is in fact left adjoint to $\mathrm{End}$. To do so, we need to do the following. \\

\noindent \textbf{1.)} We need to provide a natural transformation $\mathrm{id}_{\EMon{1}} \Rightarrow \mathrm{End} \circ (P \circ \B)$. Consider, for a given small $\infty$-category $C$ the Yoneda embedding
$$ C \hookrightarrow \Fun( C^{\op}, \Spc ). $$
This gives a natural transformation $$\mathrm{inc} \Rightarrow  P$$ of functors $\EMon{0}(\Catinf) \rightarrow \EMon{0}(\widehat{\Catinf})$. Precomposing with $\B : \EMon{1} \rightarrow \EMon{0}(\Catinf)$ and postcomposing with $\Omega : \EMon{0}(\widehat{\Catinf}) \rightarrow \EMon{1}$ produces a natural transformation of functors $\EMon{1} \rightarrow \EMon{1}$ of the form
$$ X = \Map_{\B X}( *, * ) \rightarrow \Omega( \Fun( \B X^{\op}, \Spc ), y_* ) = \mathrm{Nat}( y_*, y_* ).$$
for a given $\mathbb{E}_1$-monoid $X$. We note that by the Yoneda Lemma, this natural transformation is an equivalence. \\

\noindent \textbf{2.)} We need to show that for all $\mathbb{E}_1$-monoids $X$ and pointed presentable $\infty$-categories $(\mathcal{C},c \in \mathcal{C})$ the above unit map induces an equivalence
\[ \mathrm{Map}_{ \EMon{0}(\PrL) }( (\Fun( \B X^{\op}, \Spc ),y_*), (\mathcal{C},c) ) \xrightarrow{\sim} \mathrm{Map}_{ \EMon{1} }( X, \Omega(\mathcal{C},c) ).\]
To see this, note that $\mathrm{Map}_{ \EMon{0}(\PrL) }( \Fun( \B X^{\op}, \Spc ), (\mathcal{C},c) )$ is obtained as the fiber of
\[ \core( \Fun^L( \Fun( \B X^{\op}, \Spc ), \mathcal{C} ) ) \xrightarrow{\mathrm{ev}_*} \core( C ) \]
over $c \in \core( C )$. By the universal property of the presheaf category, we have
\[\Fun^L( \Fun( \B X^{\op}, \Spc ), \mathcal{C} ) ) \simeq \Fun( \B X, \mathcal{C} ),\]
and upon taking fibers we see that 
\[ \begin{array}{rcl}
\mathrm{Map}_{ \EMon{0}(\PrL) }( \Fun( \B X^{\op}, \Spc ), (\mathcal{C},c) ) & \simeq & \mathrm{Map}_{ \EMon{0}(\widehat{\Catinf}) }( \B X, (\mathcal{C},c) ) \\ & \simeq & \mathrm{Map}_{ \EMon{1} }( X, \Omega(\mathcal{C},c) )
\end{array}
 \]
where the last equivalence follows from the adjunction $\B \dashv \Omega$.
\end{proof}

\begin{remark} \label{identificationoffunctors}
Since the above equivalence is natural, it means in particular that for a given map of $\mathbb{E}_1$-monoids $f : X \rightarrow Y$, the functor
$$ f_! : \Fun(\B X, \Spc) \rightarrow \Fun(\B Y, \Spc)$$
given by left Kan extension identifies with the change-of-base functor
$$ f_! : \mathrm{LMod}_X(\Spc) \rightarrow \mathrm{LMod}_Y(\Spc)$$
which sends a left $X$-module $M$ to the relative tensor product $Y \otimes_X M$.
\end{remark}

We can generalize from $X$-actions in $\Spc$ to any presentable $\infty$-category.

\begin{theorem}[\cite{lurieha}, Theorem 4.8.4.6] \label{generalrmod}
Let $\mathcal{C}$ be a presentably monoidal $\infty$-category, let
$\mathcal{M}$ be a right $\mathcal{C}$-module in $\PrL$, and let
$A \in \mathrm{Alg}(\mathcal{C})$ be an algebra object. Then the canonical functor
\[
\mathcal{M} \otimes_{\mathcal{C}} \mathrm{RMod}_A(\mathcal{C})
\longrightarrow
\mathrm{RMod}_A(\mathcal{M})
\]
is an equivalence of $\infty$-categories.
\end{theorem}

\begin{corollary}
Let $X$ be an $\mathbb{E}_1$-monoid, and suppose $\mathcal{D}$ is presentable. There is a canonical equivalence of $\infty$-categories
$$ \mathrm{Fun}( \B X,  \mathcal{D} ) \simeq \mathrm{LMod}_X( \mathcal{D} ).$$
\end{corollary}

\begin{proof}
This follows at once from Theorem \ref{generalrmod} from the observation that $\mathrm{Fun}( \B X,  \mathcal{D} ) \simeq \mathrm{Fun}( \B X,  \Spc ) \otimes \mathcal{D}.$ 
\end{proof}

If $A,B$ are a pair of algebra objects in a monoidal $\infty$-category $\mathcal{C}$, it makes sense to speak of $A$-$B$-bimodules, as defined in \cite[Chapter 4.3]{lurieha}. These together form an $\infty$-category
$$  {}_A\mathrm{BMod}_B(\mathcal C). $$
Lurie shows the following characterization of bimodules. Recall that $\mathrm{LMod}_A(\mathcal C)$ is canonically a right module for $\mathcal{C}$ in $\mathrm{Pr}^L$. Hence it makes sense to speak of $B$-right modules in this category.

\begin{theorem}[\cite{lurieha}, Remark 4.3.3.8]
Let $\mathcal{C}$ be a presentably monoidal $\infty$-category, and let $A,B$ be two algebra objects in $\mathcal{C}$. There exists a canonical equivalence 
\[ {}_A\mathrm{BMod}_B(\mathcal C) \simeq \mathrm{RMod}_B\bigl(\mathrm{LMod}_A(\mathcal C)\bigr). \]
\end{theorem}

\begin{corollary}
Let $X,Y$ be two $\mathbb{E}_1$-monoids. There exists a canonical equivalence of $\infty$-categories
$$  {}_X\mathrm{BMod}_Y(\Spc) \simeq \Fun( \B X \times \B Y^{\op}, \Spc )$$
\end{corollary}

\begin{proof}
Combining the previous results we see that
\[ \begin{array}{rcl}
{}_X\mathrm{BMod}_Y(\Spc) & \simeq & \mathrm{RMod}_Y\bigl(\mathrm{LMod}_X(\Spc)\bigr) \\
 & \simeq & \Fun( \B Y^{\op}, \Fun( \B X, \Spc ) )  \\
 & \simeq & \Fun( \B X \times \B Y^{\op},  \Spc )
\end{array} \]
\end{proof}
 
In the case that $A = B$ is actually the same algebra object in $\mathcal{C}$, the $\infty$-category of bimodules
\[{}_A\mathrm{BMod}_A(\mathcal C)\]
is equipped with an $\mathbb{E}_1$-monoidal structure, given on underlying objects by $(X,Y) \mapsto X \otimes_A Y$. Hence it makes sense to talk about algebra objects in the $\infty$-category of bimodules.

\begin{theorem}[\cite{lurieha} Corollary 3.4.1.7, and Section 4.4, specifically Theorem 4.4.1.28]
Let $\mathcal{C}$ be a presentably monoidal $\infty$-category, and let $A \in \mathrm{Alg}(\mathcal{C})$. Then there is a canonical equivalence
\[
\mathrm{Alg}(\mathcal{C})_{A/}
\;\simeq\;
\mathrm{Alg}\bigl({}_A\mathrm{BMod}_A(\mathcal{C})\bigr).
\]
\end{theorem}

Now we only need to specialize to the case of $\mathcal{C}$ being $\Spc$.

\begin{corollary}
Given $X \in \EMon{1}$, there is an equivalence
$$ \mathbf{act}_X : \EMon{1}_{ X / } \simeq \EAlg{1}( \Fun( \B X \times \B X^{\op}, \Spc )) $$
that sends a map $f : X \rightarrow Y$ of $\mathbb{E}_1$-monoids to the $X$-bimodule with underlying left action given by $f^*( Y^\circlearrowleft )$.
\end{corollary}

The above situation has a symmetric monoidal version as well. If $\mathcal{C}$ is a symmetric monoidal $\infty$-category, denote by $\mathrm{CAlg}(\mathcal{C})$ the $\infty$-category of $\mathbb{E}_\infty$-algebras in $\mathcal{C}$.

\begin{theorem}[\cite{lurieha} Corollary 3.4.1.7 and Corollary 4.5.1.5]
Let $\mathcal{C}$ be a presentably symmetric monoidal $\infty$-category, and let $A \in \mathrm{CAlg}(\mathcal{C})$. Then there is a canonical equivalence
\[
\mathrm{CAlg}(\mathcal{C})_{A/}
\;\simeq\;
\mathrm{CAlg}\bigl(\mathrm{LMod}_A(\mathcal{C})\bigr).
\]
\end{theorem}

\begin{corollary} \label{abstractsymmetricmonoidalstuff}
Given $X \in \EMon{\infty}$, there is an equivalence
$$ \mathbf{act}_X : \EMon{\infty}_{ X / } \simeq \EAlg{\infty}( \Fun( \B X, \Spc )). $$
\end{corollary}

As a very special case of a map of $\mathbb{E}_1$-monoids, consider $X \rightarrow X^{\gp}$. Since $\B X^{\gp} = | \B X |$ and $|\B X|$ is the localization of $\B X$ at the set of all arrows, we get the adjunction
\[
\begin{tikzcd}
  \Fun( \B X, \Spc ) \arrow[r, bend left=10, "\B X^{-1}"] &
  \Fun( | \B X |, \Spc ) \arrow[l, bend left=10, "\rotadj", swap]
\end{tikzcd}
\]
where the right adjoint is fully faithful, and identifies spaces with $| \B X |$-actions with those spaces with $\B X$-action, such that action with any $x \in X$ is invertible. We can equip $X^{\gp}$ with its $X$-action coming from the map $X \rightarrow X^{\gp}$. It is clear that the localization sends $X^\circlearrowleft$ to $(X^{\gp})^\circlearrowleft$ since both functors are the corepresented functors of the base points, hence we obtain another description of the group completion:
$$ X^{\gp} = \B X^{-1} (X^\circlearrowleft). $$

\section{The case of commutative $\mathbb{E}_n$-monoids}

In the following section, we always assume that $X$ is at least somewhat commutative, i.e.\ an $\mathbb{E}_n$-monoid for ${n \geq 2}$.  In this case homology arguments can be quite useful, due to a version of the group completion theorem going back to Quillen, which will be the centrepiece for this section. Before we get to this, some generalities on Day convolution are needed. As discussed in the previous section, the assumption that $X$ is at least an $\mathbb{E}_2$-monoid implies that the categories $\B X, \E X$ and $\G X$ are monoidal $\infty$-categories. Thus it makes sense to talk about a highly useful monoidal structure on the categories $\Fun( \B X, \Spc )$, $\Fun( \E X, \Spc )$ and so on.

\begin{proposition} \label{dayconvolution}
Suppose $(\mathcal{C},\otimes,1)$ is an $\mathbb{E}_n$-monoidal $\infty$-category for $n \geq 1$. Then there exists a unique $\mathbb{E}_n$-monoidal structure $\otimes^{Day}$ on $\Fun(C, \Spc )$ called \emph{Day convolution} such that
\begin{enumerate}
\item $\otimes^{Day}$ preserves colimits in both variables.
\item The Yoneda embedding $\mathcal{C}^{\op} \rightarrow \Fun(C, \Spc )$ is strong $\mathbb{E}_n$-monoidal.
\item The functor $\Fun( (-)^{\op} , \Spc) : \Catinf \rightarrow \text{Pr}^L$ is strong symmetric monoidal.
\item In particular, if $f : \mathcal{C} \rightarrow \mathcal{D}$ is a strong $\mathbb{E}_n$-monoidal functor between $\mathbb{E}_n$-monoidal $\infty$-categories, then the functor
$$f_! : \Fun(\mathcal{C}, \Spc) \rightarrow \Fun(\mathcal{D}, \Spc)$$
given by left Kan extension along $f$ is strong $\mathbb{E}_n$-monoidal, and
$$f^* : \Fun(\mathcal{D}, \Spc) \rightarrow \Fun(\mathcal{C}, \Spc)$$
given by precomposition with $f$ is lax $\mathbb{E}_n$-monoidal.
\end{enumerate}
\end{proposition}

\begin{proof}
Points 1 to 3 follow immediately from Lurie \cite{lurieha}, Section 4.8.1. The first part of point 4 about $f_!$ is a consequence of 3 as the strong monoidal functor $\Fun( (-)^{\op} , \Spc)$ sends maps of algebras to maps of algebras. Lastly, the functor $f^*$ is the right adjoint to $f_!$, and thus as a right adjoint to a strong monoidal functor itself lax monoidal.
\end{proof}

Now assume $X$ is an $\mathbb{E}_n$-monoid for $n \geq 2$. Then $X^\circlearrowleft : \B X \rightarrow \Spc$ is the functor corepresented by the unit of $\B X$, hence by Proposition \ref{dayconvolution} it is the unit with respect to Day convolution. Let $Y$ be a space with a left $X$-action, i.e.\ a functor $\B X \rightarrow \Spc$. Then $Y$ is automatically a module over $X^\circlearrowleft$ with respect to Day convolution, since $X^\circlearrowleft$ is the unit. Denote the composition $ p^*(Y) = Y \circ p : \E X \rightarrow \Spc$ by $Y^\rightarrow$. Since $p^*$ is lax-monoidal, it follows that $Y^\rightarrow$ is a module over $X^\rightarrow$ with respect to Day convolution. If $E$ is a spectrum, denote by $E[Y] = E \otimes \Sigma^{\infty}_+ Y$. We will write $E_p(Y)$ for the homotopy groups $\pi_p( E[Y] )$. In the particular case when $R$ is an $\mathbb{E}_\infty$-ring spectrum and $X$ is an $\mathbb{E}_n$-monoid, the spectrum $R[X]$ has the structure of an $\mathbb{E}_n$-ring spectrum. Note that if $Y$ has an $X$-action, it is immediate that $E[Y]$ is an $\mathbb{S}[X]$-module. In particular, the adjunction
\[
\begin{tikzcd}
\Fun( \B X, \Spc ) \arrow[r, bend left=10, "\mathbb{S}{[-]}"] & \mathbb{S}[X]\text{-Mod} \arrow[l, bend left=10, "\rotadj"', "\Omega^\infty"]
\end{tikzcd}
\]
exhibits $\mathbb{S}[X]\text{-Mod}$ as the stabilization of $\Fun( \B X, \Spc )$. We need some more generalities about localizations of ring spectra.

\begin{lemma}[\cite{lurieha}, section 7.2.3.]
Let $R$ be a ring spectrum and $S \subset \pi_* R$ be a multiplicative subset of homogeneous elements satisfying the \emph{left Ore condition}. This means that
\begin{enumerate}
\item For $ x \in \pi_* R$ and $s \in S$, there exist $y \in \pi_* R$ and $t \in S$ such that $tx = ys$.
\item Let $ x \in \pi_* R$. Suppose for a given $s \in S$ we have $xs = 0$. Then there exists $t \in S$ such that $tx = 0$.
\end{enumerate}
Call an $R$-module $M$ $S$-local if multiplication with $s$ on $M$ is an equivalence for all $s \in S$. Then the inclusion of the full subcategory $R$-$\text{Mod}^{S\text{-loc}}$ of all $S$-local modules into $R$-Mod has a left adjoint
$S^{-1} : R\text{-Mod} \rightarrow R\text{-Mod}^{S\text{-loc}},$
and $\pi_* S^{-1} M \cong S^{-1} \pi_* M$.
\end{lemma}

Note that the left Ore condition is automatic if $R$ is an $\mathbb{E}_2$-ring spectrum, such as $\mathbb{S}[X]$ for $X$ an $\mathbb{E}_2$-monoid. We can apply this to the ring spectrum $\mathbb{S}[X]$, with multiplicative subset $\pi_0 (X) \subset \pi_*(\mathbb{S}[X])$. It is immediate that if $M$ is a $\pi_0(X)$-local $\mathbb{S}[X]$-module, its underlying space $\Omega^\infty M$ will be $\B X$-local, i.e.\ the action with any $x \in X$ is invertible. Hence we get a commuting square of right adjoint functors
$$\xymatrix{
\Fun( \B X, \Spc )  & \Fun( | \B X |, \Spc ) \ar[l] \\
\mathbb{S}[X]\text{-Mod} \ar[u]_{\Omega^\infty} & \mathbb{S}[X]\text{-Mod}^{\pi_0 X\text{-local}} \ar[u]_{\Omega^\infty} \ar[l].
}$$

\begin{theorem} \label{homology}
Let $n \geq 2$. Suppose $E$ is a spectrum and $Y$ a space with an $X$-action. Then there is an equivalence of $\mathbb{S}[X]$-modules, 
$$(\pi_0 X)^{-1} E[Y] \xrightarrow{\simeq } E[ \colim_{\E X} Y^\rightarrow ]$$
As a consequence, there are natural isomorphisms
$$ (\pi_0 X)^{-1} E_k(Y) \xrightarrow{\cong} E_k( \colim_{\E X} Y^\rightarrow ) $$
for all $k \in \mathbb{Z}$.
\end{theorem}

\begin{proof} This argument in its original form is due to Quillen, see \cite[p.\ 5]{10.1007/BFb0080003}. The functor $X^\rightarrow : \E X \rightarrow \Spc$ is naturally an $\mathbb{E}_1$-algebra with respect to Day convolution. Therefore, the spectrum valued functor $E[Y^\rightarrow]$ is a functor with values in $\mathbb{S}[X^\rightarrow(1)] = \mathbb{S}[X]$-modules. In particular, its homotopy groups are  $\pi_0 (\mathbb{S}[X]) \cong \mathbb{Z}[ \pi_0 X ]$-modules. Since $\colim_{ \E X } E[Y^\rightarrow]$ is $\pi_0 X$-local, we have the equivalence
$$\colim_{ \E X } (\pi_0 X)^{-1} E[Y^\rightarrow] \simeq (\pi_0 X)^{-1} \colim_{ \E X } E[Y^\rightarrow] \simeq E[ \colim_{\E X} Y^\rightarrow], $$
where we have used that localization with respect to $\pi_0 X$ is a left adjoint and thus commutes with colimits.

We can now use the colimit spectral sequence in homology,
$$ E^2_{pq} = H_p( \E X; (\pi_0 X)^{-1} E_q Y) \implies E_{p+q} ( \colim_{\E X} Y^\rightarrow ). $$
The coefficient system $(\pi_0 X)^{-1} E_q Y$ sends all arrows in $\E X$ to equivalences, hence factors through $| \E X | \simeq \pt$. All terms on the $E^2$-page of this spectral sequence thus vanish for $p > 0$, hence the spectral sequence collapses and gives the claimed isomorphisms.
\end{proof}

\begin{lemma} \label{Lisgrouplike}
For $n > 1$, there is a natural isomorphism
$$ \pi_0 L(X) \cong ( \pi_0 X )^{\gp} $$
In particular, the $\mathbb{E}_{n-1}$-monoid $L(X)$ is grouplike.
\end{lemma}

\begin{proof}
In the case of $n > 1$, the discrete monoid $\pi_0 X$ is commutative. The functor $\pi_0 : \Spc \rightarrow \Set$ commutes with colimits, hence
$$ \pi_0 L(X) \cong \colim_{ \B X } \pi_0 ( X^\circlearrowleft \times X^\circlearrowleft ) \cong \colim_{ \B (\pi_0 X) } (\pi_0 X)^\circlearrowleft \times (\pi_0 X)^\circlearrowleft, $$
and thus the statement reduces to the known fact that for (discrete) commutative monoids $M$, the group completion can be obtained as a quotient of $ M \times M$.
\end{proof}

The natural $\mathbb{E}_n$-map $X \rightarrow X^{\gp}$ induces a natural transformation of $\mathbb{E}_{n-1}$-monoids,
$$ \epsilon_X : L(X) \rightarrow L(X^{\gp}) \simeq X^{\gp} $$
which provides the inverse to $\eta_X$ if $X$ is grouplike. In the following, we denote by $U_{n-1}$ the forgetful functor $\EMon{n} \rightarrow \EMon{n-1}$. We want to argue the following theorem.

\begin{theorem} \label{Liscompletion}
Let $n \geq 2$. Then the functor \emph{$L : \EMon{n} \rightarrow \EGrp{n-1}$} is naturally equivalent to the functor
$$U_{n-1} \circ (-)^{\gp} : \EMon{n} \rightarrow \EGrp{n-1}.$$
In the case $n = \infty$, the functor $L : \EMon{\infty} \rightarrow \EGrp{\infty}$ is naturally equivalent to the group completion.
\end{theorem}

\begin{remark}
The attentive reader may ask themselves why the peculiar loss of structure happens, i.e.\ if we start with an $\mathbb{E}_n$-monoid $X$, then $L(X) = |\G X|$ is a priori only an $\mathbb{E}_{n-1}$-monoid. Theorem \ref{Liscompletion} however shows that, if $n > 1$, since $L(X)$ is equivalent to the group completion of $X$, which is an $\mathbb{E}_n$-monoid, we discover an additional $\mathbb{E}_1$-monoid structure after the fact. This extra structure cannot be visible as some $\mathbb{E}_n$-monoidal structure on the $\infty$-categories $\B X$ or $\E X$. If, for example, $\B X$ were itself $\mathbb{E}_n$-monoidal, then $X$ would be an $E_{n+1}$-monoid to begin with. The extra structure becomes more apparent when comparing the two squares
$$\xymatrix{
\G X \ar[r] \ar[d] & \E X  \ar[d] & & |\G X| \ar[r] \ar[d] & |\E X| \ar[d] \\
\E X \ar[r] & \B X & & |\E X| \ar[r] & |\B X|.
}$$
The space $|\B X|$ is the classifying space $BX$. Also, $|\E X| \simeq \pt$. Theorem \ref{Liscompletion} can be reinterpreted as saying that realization sends the particular pullback square in $\EMon{n-1}(\Catinf)$ square on the left to a pullback square in $\EMon{n-1}$. This identifies $|\G X|$ with $\Omega BX$. The additional $\mathbb{E}_1$-structure is thus induced from the general fact that loop spaces come equipped with an $\mathbb{E}_1$-structure. On a moral note, the diagram
$$\begin{tikzcd}
    \G X \arrow[r, bend left, "p^+"] \arrow[r, bend right, "p^-", swap] & \E X \arrow[l, "\Delta", swap]
\end{tikzcd}$$
can be understood as the beginning of the simplicial diagram given by the \v{C}ech nerve of the functor $\E X \rightarrow \B X$. This is an internal groupoid in the category of $\mathbb{E}_{n-1}$-monoidal $\infty$-categories, such that the object of $0$-cells, $\E X$, is a contractible $\infty$-category. We will make this more precise in the upcoming section.
\end{remark}

One may be tempted to argue that Theorem \ref{Liscompletion} follows trivially from the following two observations: $LX$ is always grouplike, and the natural transformation $\eta: X \rightarrow LX$ is an equivalence for $X$ being grouplike. However, this is not sufficient, as it only shows that the functor $L$ contains $U_{n-1} \circ (-)^{\gp}$ as a retract. One also needs the statement that $L(\eta_X) : LX \rightarrow LLX$ is an equivalence. This would follow from the statement that the maps $L(\eta_X)$ and $\eta_{L(X)} : LX \rightarrow LLX$ are homotopic, as the latter is an equivalence. We will circumvent this issue with a homology argument.

\begin{proof}
Consider the map of $\mathbb{E}_n$-monoids $X \rightarrow X^{\gp}$. Since $X^{\gp}$ is also the group completion of $X$ when viewed only as $\mathbb{E}_{n-1}$-monoids, see Corollary \ref{groupcompletionforget}, and $L(X)$ is grouplike by Lemma \ref{Lisgrouplike}, there exists a commutative triangle
\[
\begin{tikzcd}
X \arrow[d] \arrow[r]   & L(X) \\
X^{\gp}  \arrow[ru, dashed] & 
\end{tikzcd}
\]
by the universal property of the $\mathbb{E}_{n-1}$-monoidal group completion. This triangle can be interpreted as a commutative triangle in $\Fun( \B X, \Spc )$ via the functor $\text{act}_X$. As discussed before, since $| \B X | = \B X [ \B X^{-1}]$ is the localization with respect to all morphisms in $\B X$, we have the localization
\[
\begin{tikzcd}
  \Fun( \B X, \Spc ) \arrow[r, bend left=10, "\B X^{-1}"] &
  \Fun( | \B X |, \Spc ) \arrow[l, bend left=10, "\rotadj", swap]
\end{tikzcd}
\]
where $\B X^{-1}$ sends $X^\circlearrowleft$ to $(X^{\gp})^\circlearrowleft$. Stabilizing leads to the commutative square of left adjoints
$$\xymatrix{
\Fun( \B X, \Spc ) \ar[r]^{\B X^{-1}} \ar[d]_{\mathbb{S}[-]} & \Fun( | \B X |, \Spc ) \ar[d]^{\mathbb{S}[-]} \\
\mathbb{S}[X]\text{-Mod} \ar[r]^{(\pi_0 X)^{-1}} & \mathbb{S}[X]\text{-Mod}^{\pi_0 X\text{-local}}.
}$$
By Theorem \ref{homology} the spectrum $\mathbb{S}[L(X)]$ is the localization of $\mathbb{S}[X]$ at the multiplicative subset $\pi_0 X$. But this means that the map $X^{\gp} \rightarrow L(X)$ becomes an equivalence after taking $\mathbb{S}[-]$. Since it is a map of $\mathbb{E}_1$-groups, it was already an equivalence to begin with, by Corollary \ref{pluscorollary}.
\end{proof}

\begin{remark}
There is an endofunctor $\iota : \G X \rightarrow \G X$ that sends $(x,y)$ to $(y,x)$. In light of Theorem \ref{Liscompletion} it seems reasonable to expect that $\iota$ becomes the homotopy inverse on $ X^{\gp} \simeq | \G X |$ after realization. This is in fact so. Consider the inclusion $\mathrm{pos} : X \rightarrow \G X$ that sends $x$ to $(x,1)$. It is clear that there is a natural transformation $\text{const}_1 \rightarrow \mathrm{pos} \cdot (\iota \mathrm{pos} )$, which is on objects given by $(1,1) \rightarrow (x,x) = (x,1) \cdot \iota(x,1)$. Since
$$ \text{Map}_{\mathbb{E}_{n-1}}( X, | \G X | ) \simeq \text{Map}_{\mathbb{E}_{n-1}}( | \G X |, | \G X | )$$
by the universal property of the group completion, the given natural transformation must extend to a homotopy between $\text{const}_1$ and $\text{id} \cdot \iota : | \G X | \rightarrow  | \G X |$. However, we need to issue an important warning here: Thomason \cite{Thomason_80} showed that there does not in general exist a natural transformation $\text{const}_1 \rightarrow \text{id} \cdot \iota$ of functors $\G X \rightarrow \G X$. This is the case for many important examples such as $\Fin^\cong$ and $\text{Proj}(R)$ for a non-trivial discrete ring.
\end{remark}

\section{The additivity theorem}

In this section we will provide an analogue to the well-known additivity theorem in $K$-theory. The diagram
$$\begin{tikzcd}
\G X \arrow[r, bend left, "p^+"] \arrow[r, bend right, "p^-", swap] & \E X \arrow[l, "\Delta", swap]
\end{tikzcd}$$
for a given $\mathbb{E}_n$-monoid $X$ can be extended to an entire simplicial diagram $\G_\bullet(X)$, which is an internal groupoid in the category of $\mathbb{E}_{n-1}$-monoidal $\infty$-categories. The $\infty$-category
$$\G_2(X) = \G X \times_{\E X} \G X$$
can be thought of as analogous to the category $\text{Ex}(\mathcal{C})$ of exact sequences of an exact category $\mathcal{C}$. Similarly to how $K$-theory splits exact sequences, we will show that the natural functor
$$\G_2(X) \rightarrow \G X \times \G X $$
becomes an equivalence after realization. Before we prove this, let us state some general facts about internal groupoids.

\begin{definition}[ See also \cite{luriehtt} definition 6.1.2.7. and prop. 6.1.2.6 ] \label{definitioninternalgroupoids}
An \emph{internal groupoid}, in an $\infty$-category $\mathcal{C}$ with finite limits is given as a simplicial object $\mathcal{G}_\bullet : \Delta^{\op} \rightarrow \mathcal{C}$ that satisfies the \emph{groupoidal Segal condition}:
\begin{itemize}
\item For all subsets $S,S'$ of $[n] = \{0,1,\hdots,n\}$ such that $S \cup S' = [n]$ and where $S \cap S' = \set{s}$ consists of a single element $s$, the following square
$$\xymatrix{
\mathcal{G}([n]) \ar[r] \ar[d] & \mathcal{G}(S') \ar[d] \\
\mathcal{G}(S) \ar[r] & \mathcal{G}(\set{s})
}$$
is a pullback, where we think of $S$ and $S'$ as being totally ordered from the order induced by $[n]$.
\end{itemize}
We call $\mathcal{G}_0$ the \emph{underlying object} of $\mathcal{G}_\bullet$. We denote the full subcategory of $\Fun(\Delta^{\op}, \mathcal{C})$ spanned by internal groupoids by $\Grpd(\mathcal{C})$. An internal groupoid such that $\mathcal{G}_0 \simeq \pt$ is called an internal group. The subcategory of $\Grpd(\mathcal{C})$ spanned by internal groups is equivalent to the $\infty$-category $\EGrp{1}(\mathcal{C})$. (See \cite[Remark 5.2.6.5.]{lurieha}.)
\end{definition}

Looking at particular instances of the pullbacks given in Definition \ref{definitioninternalgroupoids}, one sees that internal groupoids come with a composition map
$$\circ : \mathcal{G}_1 \times_{\mathcal{G}_0} \mathcal{G}_1 \xleftarrow{\sim} \mathcal{G}_2 \xrightarrow{\delta_1} \mathcal{G}_1$$
as well as an inverse map
$$\iota : \mathcal{G}_1 \rightarrow \mathcal{G}_1$$
together with a homotopy $id_{\mathcal{G}_1} \circ \iota \simeq s_0 \delta_1$.
  
Assume furthermore that geometric realizations of simplicial objects exist in $\mathcal{C}$. Then there is an adjunction
$$\begin{tikzcd}
    \Fun(\Delta^{\op},\mathcal{C}) \arrow[r, bend left, "ev_0 \rightarrow |-|"{name=F}] & \Fun(\Delta^{1},\mathcal{C}) \arrow[l, bend left, "N^{\check{C}}"{name=G}]
\end{tikzcd}$$
where the left adjoint sends a simplicial object $X_\bullet$ to the map $X_0 \rightarrow |X_\bullet|$. The right adjoint $N^{\check{C}}$ is called \emph{\v{C}ech nerve} and sends a map $X \rightarrow Y$ to the simplicial object
$$\begin{tikzcd}
\cdots & X \times_Y X \times_Y X \arrow[r, bend right] \arrow[r] \arrow[r, bend left]  & \arrow[l, bend right=10] \arrow[l, bend left=10] X \times_Y X \arrow[r, bend right] \arrow[r, bend left]  & X \arrow[l],
\end{tikzcd}$$
with object in degree $n$ given by the $n+1$-fold pullback $ X \times_Y \cdots \times_Y X$. The \v{C}ech nerve of a map is always an internal groupoid, see \cite[Proposition 6.1.2.11]{luriehtt}.

\begin{definition} Let $X$ be an $\mathbb{E}_n$-monoid. Define $\G_\bullet(X) : \Delta^{\op} \rightarrow \EMon{n-1}(\Catinf)$ as the \v{C}ech nerve of the functor $\E X \rightarrow \B X$.
\end{definition}

Concretely, $\G_k(X)$ is obtained as the iterated pullback, taken in $\EMon{n-1}(\Catinf)$,
$$\G_k(X) = \G X \times_{\E X} \cdots \times_{\E X} \G X. $$
We obtain as special cases $\G_0(X) = \E X$ and $\G_1(X) = \G X$. The first projection $\G_k(X) \rightarrow \E X$ is equivalent under straightening to the functor
$$ X^\rightarrow \times \cdots \times X^\rightarrow : \E X \rightarrow \Spc $$
given by the $k$-fold product of $X^\rightarrow$. Similar to the cases of $\E X$ and $\G X$, the objects of $\G_k(X)$ are given by $(k+1)$-tuples of objects of $X$, and arrows $(x_i)_{0\leq i \leq k} \rightarrow (y_i)_{0\leq i \leq k}$ are determined by the datum $(z, (\alpha_i)_{0\leq i \leq k})$, where $\alpha_i : z \cdot x_i \sim y_i$ are paths in $X$. When $n \geq 2$ one can see that $| \G_k(X) |$ is always an $\mathbb{E}_{n-1}$-group. This is clear as for a given object $(x_i)_{0\leq i \leq k}$ there exists the object $(z_i)_{0\leq i \leq k}$ with
$$ z_i = x_0 \cdots x_{i-1} \cdot \hat{x}_i \cdot x_{i+1} \cdots x_k, $$
where $\hat{x}_i$ means that the $i$-th entry in the product is omitted, together with a path $(1) \rightarrow (x_i) \otimes (z_i)$. As for $\G X$, we warn the reader that this path is not part of a natural transformation in general for $k \geq 1$.

\begin{theorem}[Additivity] \label{additivitytheorem} Let $X$ be an $\mathbb{E}_n$-monoid for $n \geq 2$. For each $k \geq 0$, the canonical functor
$$\G_k(X) \rightarrow \prod_{i=1}^k \G X $$
becomes an equivalence after realization.
\end{theorem}

\begin{proof} The case $k = 0$ is clear, as the realization of $\E X$ is contractible. Now let $k \geq 1$. Under straightening the statement is equivalent to the natural map
$$ \colim_{\E X} (X^\rightarrow \times \cdots \times X^\rightarrow) \rightarrow \colim_{\E X} X^\rightarrow \times \cdots \times \colim_{\E X} X^\rightarrow $$
being an equivalence, where $\E X$ acts diagonally on the left hand side. Since both sides realize to $\mathbb{E}_1$-groups, it suffices to check that it is a homology equivalence. Using Theorem \ref{homology} the map becomes
$$ \pi_0 X^{-1} H \mathbb{Z}[X \times \cdots \times X] \rightarrow \pi_0 X^{-1} H \mathbb{Z}[X] \otimes \cdots \otimes \pi_0 X^{-1} H \mathbb{Z}[X]. $$
Note that
$$\pi_0 X^{-1} H \mathbb{Z}[X \times \cdots \times X] \simeq \pi_0 X^{-1} ( H \mathbb{Z}[X] \otimes \cdots \otimes H \mathbb{Z}[X] ).$$
The statement now follows by applying the Künneth Spectral sequence and using the following lemma.
\end{proof}

In the following, a ring $R$ is called a \emph{bialgebra}, if it is a comonoid in the category of rings. We will denote the comultiplication by $\Delta : R \rightarrow R \otimes R$. Let $S$ be a central multiplicative subset of $R$. We say that $S$ \emph{respects the coalgebra structure} if $\Delta(S)$ and $S \otimes S$ are multiplicatively equivalent subsets of $R \otimes R$, by which we mean that for any $s \in S \otimes S$ there exists $s' \in S \otimes S$ such that $s \cdot s' \in \Delta S$ and similarly with $\Delta S$ and $S \otimes S$ reversed. An example of such a bialgebra and multiplicative subset is the group ring $\mathbb{Z}[\pi_0(X)]$ with the subset $\pi_0(X)$, for $X$ being an $\mathbb{E}_2$-monoid. 

\begin{lemma} Let $R$ be a bialgebra and $S$ a central multiplicative subset that respects the coalgebra structure. Let $M_i$ for $i = 1,\cdots,k$ be a finite set of $R$-modules. Then the canonical map
$$ S^{-1} ( M_1 \otimes_{\mathbb{Z}} \cdots \otimes_{\mathbb{Z}} M_k) \rightarrow S^{-1} M_1 \otimes_{\mathbb{Z}} \cdots \otimes_{\mathbb{Z}} S^{-1} M_k $$
is an isomorphism of abelian groups, where the left hand side is equipped with the diagonal $R$-module structure.
\end{lemma}

\begin{proof}
An inverse map is given by
$$ \frac{x_1}{s_1} \otimes \cdots \otimes \frac{x_k}{s_k} \mapsto \frac{1}{s}( \hat{s}_1 x_1 \otimes \cdots \hat{s}_k x_k ) $$
where $s = s_1 \cdot s_2 \cdots s_k$ and $\hat{s}_i$ is the same product but with $s_i$ omitted.
\end{proof}

\begin{remark}
We note that the proof of Theorem \ref{additivitytheorem} would simplify greatly if the $\infty$-category $\E X$ were sifted. This need not always be the case, an example is given by $X = \Fin^\cong$, where $\E X = \text{FinInj}$ is not a sifted category, see also Remark \ref{threearrowcalculuscomment}.
\end{remark}

\begin{remark} Recall that $B X = | \B X |$. The \v{C}ech nerve of the map $\pt \rightarrow B X$ provides a model for the structure of $\Omega B X$ as an internal group object in $\mathbb{E}_{n-1}$-monoids, i.e.\ an $\mathbb{E}_n$-group. There is the natural comparison map
$$ | \G_\bullet(X) | = | N^{\check{C}}( \E X \rightarrow \B X ) | \rightarrow N^{\check{C}} | \E X \rightarrow \B X | = \Omega B X.$$
Let $S,S'$ be two subset of $[m]$ such that $S \cup S' = [m]$ and $S \cap S' = \set{s}$. By Theorem \ref{additivitytheorem} we have the canonical natural transformation from the pullback square, taken in $\EMon{n-1}(\Catinf)$,
$$\xymatrix{
\G_{[m]}(X) \ar[r] \ar[d] & \G_{S'}(X) \ar[d] \\
\G_{S}(X) \ar[r] & \G_{\set{s}}(X).
}$$
to the pullback
$$\xymatrix{
\G X^m \ar[r] \ar[d] & \G X^{|S'|-1} \ar[d] \\
\G X^{|S|-1} \ar[r] & \pt.
}$$
with maps that are level-wise equivalences after realization. Since realization commutes with finite products, this implies that $| \G_\bullet(X) |$ is again an internal groupoid. Theorem \ref{Liscompletion} implies that the natural map $| \G_\bullet(X) | \rightarrow \Omega B X$ is an equivalence in degrees $0$ and $1$, therefore it is an equivalence of internal groups.
\end{remark}

\section{How to model $\Omega X^{\gp}$} \label{loops}

In the following section, we always assume that $X$ is an $\mathbb{E}_n$-monoid for $n \geq 2$. A very important question in practice is how to get some understanding of the higher homotopy groups of $X^{\gp}$. This amounts to understanding $\Omega X^{\gp} = \text{Map}_{| \G X |}( (1,1), (1,1) )$. This space of course only depends on the path component of $(1,1)$ in $| \G X |$.

\begin{lemma} \label{pathcomponentofunit}
The path component of $(1,1)$ in $| \G X |$ can be identified as spanned by the pairs $(x,y) \in \G(X)$ such that $x$ and $y$ are \emph{stably equivalent}, which means there exists $k \in  X$ such that $k \cdot x \simeq k \cdot y$ in $X$.
\end{lemma}

\begin{proof}
It is clear that $(1,1)$ is a stably equivalent pair. If $(x,y)$ is stably equivalent, choose $k \in  X$ such that $k \cdot x \simeq k \cdot y$. We have the path
$$ (1,1) \xrightarrow{k \cdot x} ( k \cdot x, k \cdot x ) \leftarrow (x,y), $$
which means that $(x,y)$ lies in the same path component as $(1,1)$. We now claim that the subcategory of $\G X$ spanned by stably equivalent pairs is closed under ingoing and outgoing arrows. Suppose $(x,y)$ is stably equivalent, and let $k \in  X$ such that $k \cdot x \simeq k \cdot y$.
\begin{itemize}
\item  If
$$ (x,y) \xrightarrow{k',\alpha,\beta} (x',y')$$
is an arrow, then it follows that
$$k \cdot x' \simeq_{k \cdot \alpha} k \cdot k' \cdot x \simeq k' \cdot k \cdot x \simeq k' \cdot k \cdot y \simeq k \cdot k' \cdot y \simeq_{k \cdot \beta} k \cdot y', $$
hence the target is stably equivalent, where commutativity of $X$ was used in the second and fourth equivalence.
\item If 
$$ (x',y') \xrightarrow{k',\alpha,\beta} (x,y)$$
is an arrow, then it follows that
$$ k \cdot k' \cdot x' \simeq_{k \cdot \alpha} k \cdot x \simeq k \cdot y \simeq_{k \cdot \beta} k \cdot k' \cdot y', $$
hence the source is stably equivalent.
\end{itemize} 
\end{proof}

It turns out, in order to compute the mapping space $\text{Map}_{| \G X |}( (1,1), (1,1) )$ the entire path component of $(1,1)$ is not needed. Consider the subcategory $\mathbf{Bin}(X)$ of $\G X$ defined as the full subcategory spanned by objects of the form $(x,x)$ for $x \in X$. We will for simplicity identify the objects of $\mathbf{Bin}(X)$ with the objects of $X$. Note that $\E X$ forms a wide, non-full subcategory of $\mathbf{Bin}(X)$. In particular, $\mathbf{Bin}(X)$ is a connected $\infty$-category.

\begin{proposition} \label{binarycomplexes}
The inclusion of $\mathbf{Bin}(X)$ into the subcategory $\G^{st.eq} X$ of $\G X$ spanned by stably equivalent pairs is cofinal. In particular there are equivalences
$$ \Omega X^{\gp} \simeq \Omega |\mathbf{Bin}(X)| \simeq \Map_{ \mathbf{Bin}(X)[ \E X^{-1}] } ( 1, 1).  $$
\end{proposition}

\begin{proof}
The monoidal functor $\Delta : \E X \rightarrow \G X$ equips $\G X$ with an $\E X$-action via multiplication from the left. Note that this functor factors through $\mathbf{Bin}(X)$, hence $\mathbf{Bin}(X)$ is stable under this action. Now let $(x,y)$ be a stably equivalent pair and denote by $\iota$ the inclusion $\mathbf{Bin}(X)$ into $\G^{st.eq} X$. Chose $k$ such that $k \cdot x \simeq k \cdot y$ and denote the resulting arrow $(x,y) \rightarrow (k\cdot x, k \cdot x)$ by $\hat{k}$. The arrow $ 1 \rightarrow k$ in $\E X$ induces a natural transformation of functors $id_{\G X} \Rightarrow k \otimes-$. This altogether gives maps
\[
\begin{tikzcd}
(x,y) / \iota \arrow[r, bend left=10, "k \otimes -"] & (k \cdot x, k \cdot x) / \iota \arrow[l, bend left=10, "\hat{k}^*"]
\end{tikzcd}
\]
together with a homotopy from $\id_{(x,y) / \iota}$ to $\hat{k}^* \circ k \otimes -$. This implies that $| (x,y) / \iota |$ is a retract of $| (k \cdot x, k \cdot x) / \iota |$, but the latter is contractible, as $(k \cdot x, k \cdot x)$ lies in the subcategory given by $\mathbf{Bin}(X)$. Therefore $(x,y) / \iota$ is contractible for all pairs $(x,y)$ that are stably equivalent, which implies that $\iota$ is cofinal by Quillen's Theorem A.

The claimed equivalences come from the following observation. The $\mathbb{E}_\infty$-monoid $\Omega X^{\gp} = \Omega | \G X |$ is computed as $\Map_{| \G X |}( (1,1),(1,1) )$. This space only depends on the subcategory $\G^{st.eq} X$. Since $\mathbf{Bin}(X)$ sits cofinally in this $\infty$-category, we get 
$$\Map_{| \G X |}( (1,1),(1,1) ) \simeq \Map_{| \mathbf{Bin}(X) |}( 1, 1 ).$$
In order to compute $|\mathbf{Bin}(X)| = \mathbf{Bin}(X)[ \mathbf{Bin}(X)^{-1} ]$, one does not need to formally invert all morphisms. Consider the subcategory $\E X$ in $\mathbf{Bin}(X)$. Since every morphism $ x \xrightarrow{k, \alpha, \beta} y $ factors as
$$ x \xrightarrow{k, \text{id}, \text{id}} k \cdot x \xrightarrow{1, \alpha, \beta} y $$
with the first arrow in $\E X$ and the second arrow invertible, it suffices to localize by $\E X$.
\end{proof}

\begin{remark}
The name $\mathbf{Bin}(X)$ was chosen in analogy with Grayson's binary complexes \cite{Grayson12} to be an unstable/non-exact version. We recall that the objects of $\mathbf{Bin}(X)$ are given by the objects of $X$. A morphism $x \rightarrow y$ is given by the datum $(k, \alpha, \beta)$ where $k \in X$ and $\alpha : k \cdot x \rightarrow y $ and $\beta : k \cdot x \rightarrow y$ are two unrelated equivalences in $X$. The mapping spaces out of $1$ can be directly computed using Lemma \ref{mappingspaces} as
$$ \Map_{\mathbf{Bin}(X)} ( 1, x ) \cong \text{Aut}(x),$$
where $\text{Aut}(x) = \Map_X(x,x)$. It is a pleasant and informative exercise to compute $\pi_1 X^{\gp}$ from this description in the case $X = \Fin^\cong$ or $\text{Proj}(R)$ in the case of a discrete ring to recover the known groups $\pi_1 \mathbb{S}$ or $K_1(R)$. It might also be helpful to try to understand how a zig-zag in $\mathbf{Bin}(\text{Proj}(R))$ starting and ending at $0$ produces a binary complex of projective $R$-modules.
\end{remark}

\begin{remark} \label{threearrowcalculuscomment}
Warning! While it may seem as if the localization of $\mathbf{Bin}(X)$ with respect to $\E X$ should satisfy some notion of a three arrow calculus of fractions to allow the computation of $\Map_{ \mathbf{Bin}(X)[ \E X^{-1}] } ( 1, 1)$, this cannot be the case in general. Consider the example of $X = \Fin^\cong$. Recall that $\E \Fin^\cong \simeq \FinInj$. Assume that we have an equivalence
$$\Map_{ \mathbf{Bin}(\Fin^\cong)[  \FinInj^{-1}] } ( \emptyset , \emptyset) \simeq \colim_{ N \in \FinInj } \text{Aut}(N).$$
A computation reveals that $\pi_0$ of the latter space is given by the set of cycle types of permutations, which classify all permutations up to conjugacy. This set is too big. One would expect further quotient operations to occur so that e.g.\ $(123)$ is mapped to zero, giving the claimed isomorphism with $\pi_1 \mathbb{S} = \left\lbrace \pm 1 \right\rbrace$. The fundamental reason for this is the existence of switching 2-cells. The square
$$\xymatrix{
\emptyset \ar[r] \ar[d] & N \ar[d] \\
N \ar[r] & N \oplus N 
}$$
one obtains by extending the canonical map $\emptyset \xrightarrow{N, \text{id},\text{id}} N $ in $\E \Fin^\cong $ along the action of $\E \Fin^\cong $ on itself only commutes up to a choice of switching isomorphism $\tau : N \oplus N \cong N \oplus N $. This problem will disappear when one considers the colimit instead in $\mathbb{E}_1$-groups, as we will see next. Note that this difference also shows that $\E \Fin^\cong = \FinInj$ is not sifted (neither when considered as a 1-category nor as an $\infty$-category). \end{remark}

The $\infty$-category $\mathbf{Bin}(X)$ is not just a convenient simplification when computing $\Omega X^{\gp}$. Rather, the projection $p^- : \mathbf{Bin}(X) \rightarrow \E X$ is a left fibration and classifies an important functor. To see how, recall that $X^\rightarrow : \E X \rightarrow \Spc$ is an $\mathbb{E}_1$-algebra with respect to Day convolution. In particular, it has a unit map $1_{\E X} \rightarrow X^\rightarrow$. Since $1_{\E X}$ is corepresented by $1 \in \E X$ which is initial, this is equivalent to the constant functor with value a point. The natural transformation $u : \pt \rightarrow X^\rightarrow$ for given $x \in \E X$ is simply the map that sends $\pt$ to $x \in X$. Now take the pullback
$$\xymatrix{
\text{Aut}(u) \ar[r] \ar[d] & \pt \ar[d]^u \\
\pt \ar[r]^u & X^\rightarrow
}$$
in $\Fun(\E X,\Spc)$, which gives the functor $\text{Aut}(u) : \E X \rightarrow \EGrp{1}$, which sends $x \in \E X$ to $\text{Aut}(x) = \Map_X(x,x)$. Taking $B \text{Aut}(u)$ gives a subfunctor of $X^\rightarrow$, i.e.\ one that is value-wise the inclusion of the path component of $x$ in $X^\rightarrow$. It is immediate upon inspection that
$$\text{Un}(B \text{Aut}(u) \rightarrow X^\rightarrow ) \simeq (\mathbf{Bin}(X) \rightarrow \G X ),$$
hence  $\mathbf{Bin}(X) \xrightarrow{p^-} \E X$ classifies the functor $B \text{Aut}(u) : \E X \rightarrow \Spc.$

\begin{theorem}
Let $X$ be an $\mathbb{E}_n$-monoid for $n \geq 2$. There is a natural equivalence of $\mathbb{E}_1$-groups
\emph{$$ \Omega X^{\gp} \simeq \colim_{x \in \E X } \text{Aut}(x) $$}
where the colimit is taken in the $\infty$-category of $\mathbb{E}_1$-groups.
\end{theorem}

\begin{proof}
Since $\mathbf{Bin}(X) \rightarrow \E X$ classifies the functor $x \mapsto B\text{Aut}(x)$, we have the equivalence
$$ | \mathbf{Bin}(X) | \simeq \colim_{x \in \E X } B \text{Aut}(x).$$
The claim now follows from the equivalence of the $\infty$-categories of $\mathbb{E}_1$-groups and connected spaces provided by the functors $B$ and $\Omega$, together with Proposition \ref{binarycomplexes}.
\end{proof}

\begin{example}
The functor $\text{Aut}(u) : \E X \rightarrow \EGrp{1}$ can be understood on morphisms as follows. The canonical morphisms $x \xrightarrow{k, \text{id}} k\cdot x$ are mapped to the stabilization maps
\[
\begin{array}{rcl}
x  & \mapsto & \text{Aut}(x) \\
(k, \text{id}) \downarrow  & \mapsto & \hspace{2.5ex}\downarrow \id_k \cdot( - ) \\
k\cdot x  & \mapsto & \text{Aut}(k \cdot x)
\end{array}
\]
and the morphisms given by automorphisms $x \xrightarrow{1, \alpha} x$ are mapped to conjugation
\[
\begin{array}{rcl}
x  & \mapsto & \text{Aut}(x) \\
(1, \alpha) \downarrow & \mapsto & \hspace{2.5ex}\downarrow \alpha ( - ) \alpha^{-1} \\
x  & \mapsto & \text{Aut}(x).
\end{array}
\]

Since $\pi_0 \Omega X^{\gp}$ is abelian, we can compute $\pi_0$ of the colimit in question by taking the abelianization of all respective groups, i.e.\
$$ \pi_1 X^{\gp} = \pi_0 \Omega X^{\gp} \cong \pi_0 \colim_{x \in \E X } \text{Aut}(x) \cong \colim_{x \in \E X } (\pi_0 \text{Aut}(x))_{ab}. $$
But this means that the action of an arrow of type $x \xrightarrow{1, \alpha} x$ becomes trivial---Hence only the formal maps $x \xrightarrow{k, \text{id}} k\cdot x$ matter. From this we can now read off the classical formulas for the examples of
$$ \pi_1 \mathbb{S} \cong \pi_1 ( (\Fin^\cong)^\gp ) \cong (\colim_n \Sigma_n)_{ab} $$
and 
$$ K_1 R \cong \pi_1 ( \text{Proj}(R)^\gp ) \cong (\colim_n \mathrm{Gl}_n(R))_{ab} $$
for $R$ a discrete ring.
\end{example}

\section{Comparison with the telescope construction} \label{telescopes}

A more classical approach to construct the group completion of a homotopy commutative topological monoid $M$ is given as follows. Construct a space $M_\infty$ from some iterated telescope involving left multiplication by elements $m : M \rightarrow M$, where $m$ goes over a set of generators of $\pi_0 M$. This space now has the property that $\pi_0 M_\infty = (\pi_0 M)^{\gp}$, but will not be a topological monoid anymore. Then apply the plus construction to $M_\infty$ to get a model for the group completion. (See e.g.\ McDuff, Segal \cite{McDuff1976}, Randal-Williams \cite{10.1093/qmath/hat024} and Nikolaus \cite{Nikolaus_2017})

For the case of $X$ being an $\mathbb{E}_2$-monoid, the construction of $X^{\gp}$ as $\colim_{\E X} X^\rightarrow$ does not have any of these defects. In particular, the colimit will naturally be an $\mathbb{E}_2$-group and no plus construction is needed. The necessary  trade-off however is that $\E X$ will not be filtered, and might be an actual $\infty$-category in general.

Let us now assume the following. Let $D$ be a filtered $\infty$-category and $\iota : D \rightarrow \E X$ be a functor. Define $X_\infty = \colim_D X^\rightarrow \circ \iota$. Then there is a natural comparison map
$$ X_\infty = \colim_D X^\rightarrow \circ \iota \rightarrow \colim_{\E X} X^\rightarrow = X^{\gp} $$
This map need not be an equivalence. An example where this will not be the case, regardless of the choice of such a functor, is $X = \Fin^{\cong}$. Nikolaus gives several different equivalent conditions however for when the comparison map is in fact an equivalence.

\begin{proposition}[See \cite{Nikolaus_2017}, Proposition 6] \label{cyclic_triviality}
The following are equivalent:
\begin{enumerate}
\item $X_\infty$ is $\pi_0(X)$-local, that is $\pi_0(X)$ acts invertibly on $X_\infty$.
\item The map $X \rightarrow X_\infty$ exhibits $X_\infty$ as the universal $\pi_0(X)$-local space.
\item The canonical map $X_\infty \rightarrow X^{\gp}$ is an equivalence.
\item The fundamental groups of all components of $X_\infty$ are abelian.
\item The fundamental groups of all components of $X_\infty$ are hypoabelian.
\item For every $x \in X$ consider the induced map $\Fin^{\cong} \rightarrow X$ (using that $\Fin^\cong$ is the free $\mathbb{E}_\infty$-monoid on a single generator). The map
$$ \Sigma_3 \rightarrow \pi_1( X, x^3) \rightarrow \pi_1( X_\infty, x^3 )$$
has $(123)$ in its kernel. 
\item For every $x$ there is an $n \geq 2$ such that the map
$$ \Sigma_n \rightarrow \pi_1( X, x^n) \rightarrow \pi_1( X_\infty, x^n )$$
has the permutation $(123\hdots n)$ in its kernel.
\end{enumerate}
\end{proposition}

We will add the following condition to the above list. Recall that a functor between $\infty$-categories is called cofinal if restriction along it does not change colimits.

\begin{proposition} \label{cofinal_inclusion} Let $D$ be filtered. Then a functor $\iota : D \rightarrow \E X$ is cofinal iff any of the equivalent conditions in Proposition \ref{cyclic_triviality} hold.
\end{proposition}

\begin{remark}
The existence of a cofinal functor $\iota : D \rightarrow \E X$ with $D$ filtered in particular implies that $\E X$ is filtered.
\end{remark}

\begin{proof} It is clear that $\iota$ being cofinal implies $(3)$. For the converse, fix $x \in X$ and consider the pullback of functors $D \rightarrow \Spc$,
$$\xymatrix{
\Map_{\E X}(x, \iota(-) ) \ar[r] \ar[d] & X^\rightarrow \circ \iota \ar[d]^{\cdot x} \\
\pt \ar[r] & X^\rightarrow \circ \iota.
}$$
Since $D$ is filtered, taking colimits commutes with pullbacks and we get the pullback
$$\xymatrix{
\colim_D \Map_{\E X}(x, \iota(-) ) \ar[r] \ar[d] & X_\infty \ar[d]^{\cdot x} \\
\pt \ar[r] & X_\infty.
}$$
Using the colimit characterization of Quillen's Theorem A, \ref{quillensA}, we see that the map $\iota$ is cofinal iff $\cdot x : X_\infty \rightarrow X_\infty$ is an equivalence for all $x \in X$.
\end{proof}

We would like to end our discussion with how one might go about finding such a functor $D \rightarrow \E X$ with filtered $D$. For simplicity, assume $X$ is an $\mathbb{E}_\infty$-monoid. Given a generating set $S$ of $\pi_0(X)$, we get an induced map of $\mathbb{E}_\infty$-monoids $F(S) \rightarrow X$, where $F(S)$ is the free $\mathbb{E}_\infty$-monoid generated by set $S$.

\begin{lemma}
The free $\mathbb{E}_\infty$-monoid generated by a set $S$ is given by the symmetric monoidal groupoid $(\Fin/S)^\cong$, with the symmetric monoidal structure induced by the coproduct of $\Fin/S$. The action $\infty$-category $\E ( F(S))$ is equivalent to $(\Fin/S)^{\text{mono}}$, the category of finite sets over $S$ with morphisms being monomorphisms.
\end{lemma}

\begin{proof}
The algebraic theory of $\mathbb{E}_\infty$-monoids is modelled by the $(2,1)$-category $\text{Span}(\Fin)$, see \cite{cranch2010algebraic}. Hence the free $\mathbb{E}_\infty$-monoid on a finite set $N$ is given by
$$\text{Hom}_{\text{Span}(\Fin)}( 1, N ) \cong (\Fin/N)^\cong.$$
The formula for a general set $S$ follows from writing $S$ as a filtered colimit over finite sets. The formula for $\E ( F(S)) = (\Fin/S)^{\text{mono}}$ can be obtained by computing the mapping spaces using Lemma \ref{mappingspaces}.
\end{proof}

Now take $\pt \in \Fin^\cong$. It induces a map of $\mathbb{E}_1$-monoids $(\mathbb{N},+) \rightarrow \Fin^\cong$ since $\mathbb{N}$ is the free $\mathbb{E}_1$-monoid on a single generator, and hence a (non-monoidal!) functor $\phi : (\mathbb{N},\leq) \rightarrow \E (\Fin^\cong ) \simeq \text{FinInj}$. Fix this functor once and for all. Define $D(S)$ as the pullback in $\Catinf$,
$$\xymatrix{
D(S) \ar[r]^{\tilde{\phi}} \ar[d] & (\Fin/S)^{\text{mono}} \ar[d] \\
\prod_{S} ( \mathbb{N}, \leq ) \ar[r]_{\prod_S \phi} & \prod_{S} \text{FinInj}.
}$$
The category $D(S)$ is a filtered poset. If we consider $\iota : D(S) \rightarrow (\Fin/S)^{\text{mono}} \rightarrow \E X$ we now get a model for the telescope
$$X_\infty = \colim_{D(S)} X^\rightarrow $$
that Randal-Williams and Nikolaus construct. Since $D(S)$ is filtered, and $S$ is a set of generators for $\pi_0(X)$, it is easy to check that the homology of $X_\infty$ is $\pi_0(X)$-local. This means the comparison map
$$ X_\infty = \colim_D X^\rightarrow \circ \iota \rightarrow \colim_{\E X} X^\rightarrow = X^{\gp} $$
is acyclic, with target an $\mathbb{E}_\infty$-group, and hence realizes the plus construction, by Theorem \ref{plusconstruction}.

The case of an $\mathbb{E}_2$-monoid can be done similarly, but one needs to work with free $\mathbb{E}_2$-monoids instead of free $\mathbb{E}_\infty$-monoids, and choose a well-ordering on the set $S$. This amounts to replacing the category $\Fin$ by the braided monoidal category of vines, see \cite{doi:10.1080/00927879708825919}.

\begin{example}
Suppose $X$ is an $\mathbb{E}_\infty$-monoid such that for all $x \in X$ there exists a path from the cyclic permutation $(123)$ of $x \cdot x \cdot x$ to the identity of $x \cdot x \cdot x$ in $\text{Aut}(X)$. Then the approximation of $\E X$ via $D(S)$ satisfies condition $(6)$ in Proposition \ref{cyclic_triviality}. It follows that $D(S) \rightarrow \E X$ is cofinal, from which we conclude that $\E X$ is filtered. This is in particular true if $A$ is a real or complex Banach algebra, and $X = \text{Proj}^{top}(A)$ is the $\infty$-category obtained by taking finitely generated projective $A$-modules and equipping the mapping spaces with their usual topology. In the complex case, let $U_\infty (A)$ be the colimit of $U_n(A) = \text{Aut}( A^n )$. It follows from the results of Section \ref{loops} together with cofinality of the inclusion $\text{Free}^{top}(A) \rightarrow \text{Proj}^{top}(A)$ of finitely generated free $A$-modules, see Theorem \ref{cofinalitytheorem}, that
$$ K^{top}(A)_{\geq 0} = \text{Proj}^{top}(A)^{\gp} \cong K_0(A) \times B ( \colim_{ P \in \E X } \text{Aut}( P )) \cong K_0(A) \times B U_\infty (A), $$
and similarly with $O_\infty (A)$ and $O_n(A)$ for the real case.
\end{example}

\section{A relative version of the group completion} \label{relative}

Often in practice, it is useful to not only have a model for the group completion of a single object as a realization, but also a relative version for the cofiber of a map induced on group completions. This works well if we deal with $\mathbb{E}_\infty$-monoids. In fact, when working with $\mathbb{E}_\infty$-monoids, we can even model cofibers of maps of $\mathbb{E}_\infty$-monoids in the $\infty$-category of $\mathbb{E}_\infty$-monoids as realizations. Let $f : X \rightarrow Y$ be a map of $\mathbb{E}_\infty$-monoids. The main point of this section will be the claim that there are symmetric monoidal $\infty$-categories $\E (f)$ and $\G (f)$, such that there are natural cofiber sequences
$$ X \rightarrow Y \rightarrow |\E (f)|$$
and
$$ X^{\gp} \rightarrow Y^{\gp} \rightarrow |\G (f)|$$
of $\mathbb{E}_\infty$-monoids and $\mathbb{E}_\infty$-groups respectively. We will give a concrete description of these two symmetric monoidal $\infty$-categories.

\begin{definition}
Let $f : X \rightarrow Y$ be a map of $\mathbb{E}_n$-monoids. Define $\E (f)$ as the pullback
$$\xymatrix{
\E (f) \ar[r] \ar[d] & \E Y \ar[d] \\
\B X \ar[r]^{\B f} & \B Y,
}$$
taken in $\EMon{n-1}(\Catinf)$.
\end{definition}

The underlying $\infty$-category $\E (f)$ is equivalent to the unstraightening of the functor $\text{act}_X( Y ) : \B X \rightarrow \Spc$, given by the induced $X$-action on $Y$ via $f$. The objects of $\E (f)$ are given by the objects of $Y$. 

\begin{lemma} \label{mappingspaceofrelative}
The mapping spaces of $\E (f)$ are given as the pullback
$$\xymatrix{
\Map_{\E (f)}( y_1, y_2 ) \ar[r] \ar[d] & X \ar[d]^{f(-)\cdot y_1} \\
\pt \ar[r]^{y_2} & Y
}$$
in $\Spc$.
\end{lemma}

\begin{proof}
This is immediate from the mapping spaces of a pullback from using the formula in Lemma \ref{mappingspaces}.
\end{proof}

\begin{example} The following are helpful special cases to keep in mind:
\begin{itemize}
\item In the case of $f$ being the zero map $ 0 : X \rightarrow 0$, it is immediate that $\E( X \rightarrow 0 )$ has a single object, and mapping space given by $X$. Hence $\E( X \rightarrow 0 ) = \B X$.
\item In the case of $f$ being the identity map $ \text{id}_X : X \rightarrow X$, it is immediate that $\E( \text{id}_X ) = \E X$.
\item In the case of $f$ being the zero map $ 0 : 0 \rightarrow X$ we recover $\E( 0 \rightarrow X ) = X$.
\item In the case of $f$ being the diagonal map $\Delta : X \rightarrow X \times X$, we immediately see that $\E ( \Delta ) = \G X$.
\end{itemize}
\end{example}

We note that there is a canonical $\mathbb{E}_{n-1}$-monoidal functor $\mathbf{can} : Y \rightarrow \E (f)$ induced from the pullback, taken in $\EMon{n-1}(\Catinf)$,
$$\xymatrix{
Y \ar[r] \ar[d] & \E (f) \ar[d] \\
\pt \ar[r] & \B X.
}$$

\begin{theorem} \label{Efmodelscofiber}
Let $f : X \rightarrow Y$ be  a map of \emph{$\mathbb{E}_\infty$}-monoids. Then there is a natural equivalence of \emph{$\mathbb{E}_\infty$}-monoids
\emph{$$ \cof ( f ) \simeq | \E (f) | $$}
induced by the functor $\mathbf{can} : Y \rightarrow \E (f)$.
\end{theorem}

This theorem follows immediately from the following lemma, using that $E (f)$ classifies $\text{act}_X( Y )$.

\begin{lemma}
Let $X$ be an $\mathbb{E}_\infty$-monoid. The following triangle of functors commutes.
\[
\begin{tikzcd}
    {\EMon{\infty}}_{ X / } \arrow{rr}{\text{act}_X} \arrow[swap]{dr}{\cof} & & \EAlg{\infty}( \Fun(\B X, \Spc ) ) \arrow{dl}{\colim^\Spc} \\
     & \EMon{\infty}. &
\end{tikzcd}
\]
\end{lemma}

For the equivalence provided by $\text{act}_X$, see Corollary \ref{abstractsymmetricmonoidalstuff}.

\begin{proof}
Note that the cofiber $\mathrm{cof}(f)$ of a map $f : X \rightarrow Y$ of $\mathbb{E}_\infty$-monoids is by definition a pushout
$$\xymatrix{
X \ar[r]^f \ar[d] & Y \ar[d] \\
0 \ar[r] & \mathrm{cof}(f)
}$$
in $\EMon{\infty}$. Pushouts of commutative algebra objects are computed as relative tensor products (see \cite{lurieha} Proposition 3.2.4.7 and Corollary 3.4.1.7), hence we see that $\mathrm{cof}(f)$ agrees with the base change of $Y$, viewed as an $X$-module, along the map $X \rightarrow 0$. By Remark \ref{identificationoffunctors}, this base change is identified with left Kan extension along the functor $\B X \rightarrow \mathrm{pt}$, which is equivalently the colimit functor $\colim^\Spc : \Fun(\B X, \Spc ) \rightarrow \Spc$.
\end{proof}

We now want to consider the question of modelling the group completion of the cofiber of a map $f : X \rightarrow Y$ of $\mathbb{E}_\infty$-monoids as a realization. We will reduce this problem to the previous case of a cofiber of a map of $\mathbb{E}_\infty$-monoids, using the observation that for an $\mathbb{E}_\infty$-monoid $Y$ it is true that
$$Y^{\gp} \simeq \cof( Y \rightarrow Y \times Y ),$$
as explained by the following proposition. Before we do so, let us remark on a feature that distinguishes the case of $\mathbb{E}_\infty$-monoids from that of $\mathbb{E}_n$-monoids: \emph{Preadditivity.}

\begin{definition}
An $\infty$-category $\mathcal C$ is called \emph{pointed} if it admits a
zero object, i.e. an object $0 \in \mathcal C$ which is both initial and terminal. A pointed $\infty$-category $\mathcal C$ is called \emph{preadditive} if it admits finite coproducts and finite products, and
for every pair of objects $C_1, C_2 \in \mathcal C$, the canonical morphism
  \[
    C_1 \sqcup C_2 \longrightarrow C_1 \times C_2
  \]
is an equivalence.
\end{definition}

The $\infty$-category $\EMon{\infty}$ and its subcategory $\EGrp{\infty}$ are both semi-additive, \cite[Proposition 2.3 and Corollary 2.5]{GepnerGrothNikolaus2015}. We will use the notation $X \oplus Y = X \times Y$ to refer to the product, or equivalently coproduct, of two $\mathbb{E}_\infty$-monoids.

\begin{proposition} \label{cofibersarecofibers}
Let $X, Y$ be $\mathbb{E}_\infty$-monoids, and $f : X \rightarrow Y$ a map of $\mathbb{E}_\infty$-monoids.
\begin{itemize}
\item The map $Y \rightarrow \cof( Y \rightarrow Y \oplus Y )$ induced by taking vertical cofibers of the square
$$\xymatrix{
0 \ar[r] \ar[d] & Y \ar[d]^\Delta \\
Y \ar[r]^{i_1} & Y \oplus Y.
}$$
identifies $\cof( Y \rightarrow Y \times Y )$ with $Y^{\gp}$.
\item Let $M(f) : X \oplus Y \rightarrow Y \oplus Y$ be the map given by the matrix
$$M_f = \begin{pmatrix}
f & 1 \\
0 & 1 
\end{pmatrix}.$$
The map $\cof(f) \rightarrow \cof(M(f))$ induced by taking vertical cofibers of the square
$$\xymatrix{
X \ar[r]^{i_1} \ar[d]^f & X \oplus Y \ar[d]^{M_f} \\
Y \ar[r]^{i_1} & Y \oplus Y.
}$$
identifies $\cof(M(f))$ with $\cof(f^{\gp})$.
\end{itemize}
\end{proposition}

\begin{proof}
Group completion is a left adjoint, hence
\[
\cof\bigl(Y \xrightarrow{\Delta} Y\oplus Y\bigr)^{\gp}
\simeq
\cof\bigl(Y^{\gp}\xrightarrow{\Delta}
Y^{\gp}\oplus Y^{\gp}\bigr).
\]
The difference map
\[
(1,-1)\colon Y^{\gp}\oplus Y^{\gp}\longrightarrow Y^{\gp}
\]
vanishes on the diagonal and therefore induces an equivalence
\[
\cof\bigl(Y^{\gp}\xrightarrow{\Delta}
Y^{\gp}\oplus Y^{\gp}\bigr)
\xrightarrow{\simeq}Y^{\gp}.
\]
Its inverse is induced by the inclusion of the first summand.

Moreover, $\cof(Y\xrightarrow{\Delta}Y\oplus Y)$ is already
grouplike. Indeed, on $\pi_0$ the class of $(a,b)$ has inverse the
class of $(b,a)$, since their sum $(a+b,a+b)$ lies in the image of
the diagonal. Thus
\[
\cof\bigl(Y\xrightarrow{\Delta}Y\oplus Y\bigr)
\simeq Y^{\gp}.
\]
Under this equivalence, the map induced by the inclusion of the
first summand agrees with the group-completion map
$Y\rightarrow Y^{\gp}$.

For the second claim, observe that we have a map of cofiber sequences of $\mathbb{E}_\infty$-monoids
$$\xymatrix{
Y \ar[r]^{i_2} \ar[d]^= & X \oplus Y \ar[d]^{M(f)} \ar[r]^{p_1} & X \ar[d]^{f} \\
Y \ar[r]^{\Delta} & Y \oplus Y \ar[r] & Y^{\gp},
}$$
hence by taking vertical cofibers we get the equivalence $\cof(M(f)) \simeq \cof(f : X \rightarrow Y^{\gp})$. It remains to show that the comparison map
$$ \cof(f : X \rightarrow Y^{\gp}) \rightarrow \cof(f : X^{\gp} \rightarrow Y^{\gp})$$
is an equivalence. Again, using that group completion commutes with colimits, it is clear that this map is an equivalence after group completion. But the left hand side is already group complete, as can be seen by noting that $\pi_0 : \EMon{\infty} \rightarrow \EMon{\infty}(Set)$ is a left adjoint, and hence preserves colimits, and the fact that any quotient of the abelian group $\pi_0(Y^{\gp})$ in the category of discrete commutative monoids remains an abelian group.
\end{proof}

We are now ready to give the following definition.

\begin{definition}
Let $f : X \rightarrow Y$ be a map of $\mathbb{E}_n$-monoids. Define $\G (f)$ as the pullback
$$\xymatrix{
\G (f) \ar[r] \ar[d] & \E ( Y \times Y ) \ar[d] \\
\B ( X \times Y ) \ar[r]^{\B (M_f)} & \B ( Y \times Y ),
}$$
taken in $\EMon{n-1}(\Catinf)$.
\end{definition}

Note that $\G (f) = \E (M_f)$. The inclusion $Y \rightarrow X \times Y$ induces the pullback
$$\xymatrix{
\G Y \ar[r] \ar[d] & \G (f) \ar[d] \\
\B Y \ar[r] & \B ( X \times Y )
}$$
in $\EMon{n-1}(\Catinf)$. Hence we see that the objects of $\G (f)$ are given as the objects of $Y \times Y$. The mapping spaces can be computed in the following way, which is an immediate application of Lemma \ref{mappingspaceofrelative}.

\begin{lemma}
The mapping spaces of $\G (f)$ are given as the pullback
$$\xymatrix{
\Map_{\G (f)}( (y_1, y_2), (y_1', y_2' )) \ar[r] \ar[d] & X \times Y \ar[d]^{M_f(-) \cdot (y_1,y_2)} \\
\pt \ar[r]^{(y_1',y_2')} & Y \times Y
}$$
in $\Spc$.
\end{lemma}

The following theorem now follows immediately from Theorem \ref{Efmodelscofiber} and Lemma \ref{cofibersarecofibers}.

\begin{theorem}
Let $f : X \rightarrow Y$ be a map of $\mathbb{E}_\infty$-monoids. Then there is a natural equivalence of $\mathbb{E}_\infty$-groups
$$ \cof ( f^{\gp} ) \simeq | \G (f) | $$
induced by the functor $\G Y \rightarrow \G (f)$.
\end{theorem}

\subsection{Application: The cofinality theorem}

\begin{definition}
Let $\iota : X \rightarrow Y$ be a map of \emph{$\mathbb{E}_n$}-monoids that is a fully faithful map on underlying spaces, i.e.\ equivalent to an inclusion of a set of path components of $Y$. We call $\iota$ \emph{cofinal} if the natural comparison map
$$ \colim_{ \E X } Y^\rightarrow \circ \iota \rightarrow \colim_{\E Y} Y^\rightarrow  $$
is an equivalence.
\end{definition}

A sufficient condition for $\iota : X \rightarrow Y$ being cofinal is given by the following condition:
\begin{itemize}
\item For all $y \in Y$, there exists $y' \in Y$, and $x \in X$ together with a path $y \cdot y' \sim \iota(x)$.
\end{itemize}
This generalizes the condition given in \cite[p.\ 225]{10.1007/BFb0080003}. In this case, using Theorem \ref{homology}, we see that the natural comparison map
$$ \colim_{ \E X } Y^\rightarrow \circ \iota \rightarrow \colim_{\E Y} Y^\rightarrow  $$
is an acyclic map between $\mathbb{E}_1$-groups and hence an equivalence. This is because the action by a given $y \in Y$ on $(\pi_0 X)^{-1} \mathbb{S}[Y]$ factors from the left and from the right through the action of some $x \in X$, which is an equivalence.

\begin{theorem}[Cofinality theorem] \label{cofinalitytheorem}
Let $\iota : X \rightarrow Y$ be a cofinal inclusion of \emph{$\mathbb{E}_n$}-monoids for $n \geq 2$. Then the cofiber $\cof(\iota^{\gp})$ of
$\iota^{\gp} : X^{\gp} \rightarrow Y^{\gp}$ is a discrete
abelian group. In particular,
$$ \pi_i X^{\gp} \rightarrow \pi_i Y^{\gp} $$
is injective in degree $0$ and an isomorphism for $i > 0$.
\end{theorem}

\begin{proof}
Cofinality of the inclusion $\iota : X \rightarrow Y$ implies that
$$\colim_{ \E X } Y^\rightarrow \circ \iota \simeq \colim_{\E Y} Y^\rightarrow . $$
Recall that $\E (\iota)$ is the unstraightening of the functor $Y^\circlearrowleft \circ \iota$. It follows that the unstraightening of the functor $Y^\rightarrow \circ \iota$ is given by the pullback
$$\xymatrix{
\text{Un}( Y^\rightarrow \circ \iota ) \ar[r] \ar[d] & \E (\iota) \ar[d] \\
\E X \ar[r] & \B X
}$$
in $\EMon{n-1}(\Catinf)$. We would like to show the following claim: There is a fully faithful inclusion $\G X \rightarrow \text{Un}( Y^\rightarrow \circ \iota )$. This can be seen as follows. Stacking the pullbacks
$$\xymatrix{
\G X \ar[r] \ar[d] & \E X \ar[d] \\
\text{Un}( Y^\rightarrow \circ \iota ) \ar[r] \ar[d] & \E (\iota) \ar[d] \\
\E X \ar[r] & \B X.
}$$
produces the wanted functor $\G X \rightarrow \text{Un}( Y^\rightarrow \circ \iota )$. To show that this functor is fully faithful, it suffices to show that the natural functor $\E X \rightarrow \E (\iota)$ is fully faithful. The mapping spaces of the latter are computed for $x_1, x_2 \in X$ as 
$$\xymatrix{
\Map_{\E (\iota)}( \iota(x_1), \iota(x_2) ) \ar[r] \ar[d] & X \ar[d]^{\iota(-)\cdot \iota(x_1)} \\
\pt \ar[r]^{\iota(x_2)} & Y
}$$
But the fiber of a map only depends on the path component of the chose base point, hence this mapping space is equivalent to the mapping space of $\E X$, as $X \rightarrow Y$ is fully faithful.

We are left to show that the inclusion $\G X \rightarrow \text{Un}( Y^\rightarrow \circ \iota )$ induces an equivalence on the path component of $(1,1) \in | \G X |$. For this we identify the objects and morphisms of $\text{Un}( Y^\rightarrow \circ \iota )$. Objects are given by pairs $(x,y)$ with $x \in X$ and $y$ in $Y$, and morphisms $(x,y) \rightarrow (x',y')$ are given by the datum $(k,\alpha, \beta)$, where $k \in X$, $\alpha : k \cdot x \sim x'$ is a path in $X$, and $\beta : \iota(k)\cdot y \sim y'$ is a path in $Y$.

Similarly to Lemma \ref{pathcomponentofunit}, it is thus clear that the path component of $(1,1) \in |\text{Un}( Y^\rightarrow \circ \iota )|$ is given by the pairs $(x,y)$ that are $X$-stably equivalent, i.e.\ there exists $k \in X$ such that $\iota(k\cdot x)\sim \iota(k)\cdot y$. Analogous to Proposition \ref{binarycomplexes}, we can see that the inclusion of the full subcategory spanned by the pairs $(x,\iota(x))$ is cofinal in this path component.

Hence we get an equivalence $\Omega | \G X | \simeq \Omega | \text{Un}( Y^\rightarrow \circ \iota ) | \simeq \Omega | \G Y |$, which concludes the proof.
\end{proof}

\begin{example}
Let $R$ be a ring and $S$ a central multiplicative subset of $R$. Write $S^{-1} \text{Proj}(R)$ for the full subcategory of $\text{Proj}(S^{-1}R)$ spanned by the modules of the form $S^{-1} P$ with $P$ f.g. projective over $R$.  Then the inclusion 
$$S^{-1} \text{Proj}(R) \rightarrow \text{Proj}(S^{-1}R)$$
is cofinal (as its image will contain all free modules). Consider the map
$$ f_S : \text{Proj}(R) \rightarrow S^{-1} \text{Proj}(R).$$
Since by definition $ f_S$ is surjective on $\pi_0$, the cofiber $\cof(f_S)$ is path-connected and hence an $\mathbb{E}_\infty$-group. Putting things together, we have a fiber sequence
$$ \Omega \cof(f_S) \rightarrow K( R ) \rightarrow K( S^{-1} R). $$
The $\mathbb{E}_\infty$-group $\cof(f_S)$ is given by the realization of $\E (f_S)$. To illustrate how this matches up with the more classical fiber sequence where the left hand term is modelled by the $K$-theory of $S$-torsion modules, let $M$ be a f.g. projective $R$-module that is $S$-torsion. Then the arrow
$$ 0 \xrightarrow{ M } S^{-1} M \cong 0 $$
is a loop in $|\E (f_S)|$. More generally, if $M$ is $S$-torsion and has a finite-length resolution by projective $R$-modules, we can use these modules to describe a loop from $0$ to itself.
\end{example}

\subsection{Application: Group completion with finite coefficients}

In some applications, one considers $K$-theory with finite coefficients rather than $K$-theory itself. The associated $K$-theory space then admits a reasonably simple description as a realization, which we explain in this section.

Let $X$ be an $\mathbb{E}_\infty$-monoid, written additively. There exists an $\mathbb{E}_\infty$-endomap $n\cdot : X \rightarrow X$ that sends $x$ to $n\cdot x = x + \cdots + x$. To see that this is a map of $\mathbb{E}_\infty$-monoids, note that the algebraic theory of $\mathbb{E}_\infty$-monoids is described by $\text{Span}(\Fin)$. There is a natural transformation
$$ \text{id}_{\text{Span}(\Fin)} \xrightarrow{n} \text{id}_{\text{Span}(\Fin)} $$
given by the span
$$ M \leftarrow \bigoplus_{i = 1}^n M \rightarrow M $$
for a finite set $M$.\footnote{Since $\text{Span}(\Fin)$ is a $(2,1)$-category, it is not too difficult to check that this assignment is in fact a natural transformation. The crucial ingredient is that given a function between sets $L \rightarrow K$, there is a canonical equivalence $$\bigoplus_{i = 1}^n L \cong \left(\bigoplus_{i = 1}^n K\right) \times_K L.$$ }

If $\mathcal{C}$ is an $\infty$-category with finite products, $\mathbb{E}_\infty$-monoids in $\mathcal{C}$ are described as finite product preserving functors $X : \text{Span}(\Fin) \rightarrow \mathcal{C}$. It is then clear that precomposition with the natural transformation given above induces the natural map $n\cdot : X \rightarrow X$.

\begin{lemma} Let $X$ be an $\mathbb{E}_\infty$-monoid. Then the cofiber $\text{cof}(n\cdot )$ of $n\cdot : X \rightarrow X$ is group complete.
\end{lemma}

\begin{proof}
Since $\text{cof}(n\cdot )$ is group complete iff $\pi_0 \text{cof}(n\cdot ) = \text{cok}(n\cdot : \pi_0 X \rightarrow \pi_0 X )$ is group complete, it suffices to show the claim for discrete commutative monoids. But on the cokernel it holds that
$$ 0 = n \cdot x = (n-1) \cdot x + x, $$
hence every element has an inverse.
\end{proof}

We denote by $( X ; \mathbb{Z}/n )^{\gp}$ the cofiber of $n \cdot : X^{\gp} \rightarrow X^{\gp}$, which by the above lemma is equivalently given by $\text{cof}( n\cdot : X \rightarrow X)$ . The following is immediate from Theorem \ref{Efmodelscofiber}.

\begin{corollary}
There is a natural equivalence 
$$( X ; \mathbb{Z}/n )^{\gp} \simeq | \E ( n \cdot : X \rightarrow X ) |.$$
\end{corollary}

\begin{example}
Suppose $S$ is a symmetric monoidal groupoid that satisfies the conditions that
  \begin{equation}
    \forall x, s \in S: \text{The translation action } s+: \text{Aut}(x) \rightarrow \text{Aut}(s+x) \text{ is injective,}
  \end{equation}
    as well as
  \begin{equation}
    \forall x \in S: \text{The multiplication action } n \cdot: \text{Aut}(x) \rightarrow \text{Aut}(n \cdot x) \text{ is injective.}
  \end{equation}
It is a straightforward application of Lemma \ref{mappingspaceofrelative} that $\E ( n \cdot : S \rightarrow S )$ is a symmetric monoidal 1-category. This is the case for example for $S = \text{Proj}(R)^\cong$ for a discrete ring $R$. Hence we get the result that the $K$-theory space with finite coefficients
$$ K( R; \mathbb{Z}/n ) \cong | \E ( n \cdot : \text{Proj}(R)^\cong \rightarrow \text{Proj}(R)^\cong ) | $$
can be described concretely as the realization of a $1$-category.
\end{example}

\begingroup
\setlength{\emergencystretch}{8em}
\printbibliography
\endgroup

\end{document}